\documentstyle[12pt,amsfonts,ifthen]{article}

\newtheorem{theorem}{Theorem}[section]
\newtheorem{lemma}[theorem]{Lemma}
\newtheorem{corollary}[theorem]{Corollary}
\newtheorem{proposition}[theorem]{Proposition}

\newtheorem{definition}[theorem]{Definition}
\newcommand{\bd}[1]{\begin{definition}\label{#1}\rm}
\newcommand{\ed}{\end{definition}}
\newcommand{\bt}[1]{\begin{theorem}\label{#1}}
\newcommand{\et}{\end{theorem}}
\newcommand{\bprop}[1]{\begin{proposition}\label{#1}}
\newcommand{\eprop}{\end{proposition}}
\newcommand{\bcor}[1]{\begin{corollary}\label{#1}}
\newcommand{\ecor}{\end{corollary}}

\newcommand{\D}{\displaystyle}

\newcommand{\lra}{\longrightarrow}
\newcommand{\Ra}{\Longrightarrow}
\newcommand{\stack}[2]{\raisebox{-2pt}{\renewcommand{\arraystretch}{.01}                                               
\begin{tabular}{c} 
$#2$\\ $ \scriptscriptstyle #1$ \end{tabular}   }}

\newcommand{\vp}{\varphi}
\newcommand{\ve}{\varepsilon} 
\newcommand{\sbnu}{\mbox{\scriptsize\boldmath${\nu}$}} 
\newcommand{\nid}{\noindent}
\newcommand{\qed}{\hfill$\Box$}
\def\1{\, {\rm I}\mskip-10mu 1}
\def\g{\, \gamma\mskip-11.8mu \gamma}\, 

\def\Del{\, \Delta{\hspace{-3mm}\Delta}}
 
\def\nn{|||} 

\def\bnu{\mbox{\boldmath${\nu}$}}

\renewcommand{\t}[1]{\tilde{#1}}

\begin{document}
\title{{Weak Convergence of $n$-Particle Systems Using Bilinear Forms}}
\par
\author{J\"org-Uwe L\"obus 
\\ Matematiska institutionen \\ 
Link\"opings universitet \\ 
SE-581 83 Link\"oping \\ 
Sverige 
}
\date{}
\maketitle
{\footnotesize
\noindent
\begin{quote}
{\bf Abstract}
The paper is concerned with the weak convergence of 
$n$-particle processes to deterministic stationary paths as 
$n\to\infty$. A Mosco type convergence of a class of bilinear 
forms is introduced. The Mosco type convergence of bilinear 
forms results in a certain convergence of the resolvents of the 
$n$-particle systems. Based on this convergence a criterion 
in order to verify weak convergence of invariant measures is 
established. Under additional conditions weak convergence of 
stationary $n$-particle processes to stationary deterministic 
paths is proved. The method is applied to the particle 
approximation of a Ginzburg-Landau type diffusion. 

The present paper is in close relation to the paper \cite{Lo11}. 
Different definitions of bilinear forms and versions of Mosco 
type convergence are introduced. Both papers demonstrate that 
the choice of the form and the type of convergence relates to 
the particular particle system. 
\bigskip

\noindent
{\bf AMS subject classification (2000)} primary 47D07, secondary 
60K35, 60J35 

\noindent
{\bf Keywords} 
Bilinear forms, convergence, Ginzburg-Landau type diffusion 
\end{quote}
}

\section{Introduction}
\setcounter{equation}{0}

The aim of the paper is to introduce a method in order to prove weak 
convergence of $n$-particle processes to deterministic stationary paths 
as $n\to\infty$. More precisely, we deal with empirical measure valued 
stochastic processes which describe the dynamics of particle configurations  
of size $n\in {\Bbb N}$. For $n\to\infty$ we establish a weak limit of 
their probability laws on the Skorohod space $D_E[0,\infty)$ where $E$ 
is a suitable space of probability measures. The limit is the degenerate 
distribution on a single probability measure valued trajectory constant 
in time. 

It is supposed that 
\begin{itemize} 
\item[{-}] the $n$-particle systems considered as measure valued stochastic 
processes are Markov with invariant probability measures $\bnu_n$, 
\item[{-}] the measure valued empirical $n$-particle processes are 
associated with strongly continuous semigroups of continuous operators 
on the $L^2$-spaces relative to the measures $\bnu_n$, 
\item[{-}] the evolution of such an $n$-particle system is initiated 
with a distribution $\psi_n\bnu_n$ over the configurations of the 
$n$ particles where the densities $\psi_n$ are uniformly bounded in 
$n\in {\Bbb N}$, cf. Theorem \ref{Theorem3.5} and Corollary 
\ref{Corollary3.6}, 
\item[{-}] capacities of sets $\{|g-\beta G_{n,\beta}g|\ge\ve\}$ 
($G_{n,\beta}$ denoting the resolvent of the $n$-particle process) 
are majorized by terms of the form $\xi (\langle g-\beta G_{n, 
\beta}g\, , \, g-\beta G_{n,\beta}g\rangle_{n}^{1/2}\, , g \, , n)$ 
which tend to zero as $n\to\infty$, cf. condition (${\cal C}4$) in 
Section 3 and Theorem \ref{Theorem3.5}.  
\end{itemize} 
The emphasis of the paper is on the method. Its usefulness is 
demonstrated by verifying weak convergence to a stationary path for 
a particle system approximating a Ginzburg-Landau type diffusion. 
\medskip

The paper is organized as follows: In Section 2, a Mosco type 
convergence of a class of bilinear forms is established. This Mosco 
type convergence of bilinear forms results in a certain convergence 
of resolvents. 

In Section 3, we relate the result of Section 2 to the class of 
particle systems the paper is concerned with. In particular, in 
Subsection 3.3, we are interested in weak convergence of invariant 
measures. Then in Subsection 3.4, we prepare the material in order 
to use the Kurtz criterion to prove relative compactness of the 
$n$-particle processes (Theorem \ref{Theorem3.5}). 

Section 4 is finally devoted to the discussion of the above mentioned 
example. 
\medskip 

The setting of Section 2 is rather general. It is more general 
than what is used in order to establish relative compactness and weak 
convergence of particle systems in Subsection 3.4 and in Section 4. 
It prepares the calculus developed in Section 2 of \cite{Lo11} and 
it is appropriate in order to show weak convergence of invariant 
measures in Subsections 3.3 and 4.1 of the present paper. 

It is also designed in order 
to establish an analysis on sequences of $L^2$-spaces and to prove 
convergence of processes on a sequence $E_n$ of state spaces. In 
particular, we mention that it is in the nature of the present work 
to be restricted to limits being stationary non-random measure valued 
paths. In fact, up to the perturbation by an $n$-particle initial 
density $\psi_n$, we are dealing with sequences of stationary 
empirical $n$-particle processes converging as $n\to\infty$ to the 
stationary solution of a partial differential equation. 
\medskip

The results of Section 2 and Subsection 3.1 might be of independent 
interest. Relations to a recent theory presented in K. Kuwae, T. Shioya 
\cite{KS03} are discussed in Subsection 3.2. We also would like to refer 
to A. V. Kolesnikov's work \cite{Ko05} and \cite{Ko06}. These papers 
develop the approach to convergence in sequences of Hilbert spaces in 
the sense of K. Kuwae and T. Shioya \cite{KS03}. For an earlier adaption 
of Mosco type convergence to non-symmetric Dirichlet forms we would like 
to draw the readers attention to the paper by M. Hino \cite{Hi98}.

\section{Convergence of Bilinear Forms}
\setcounter{equation}{0}

In this section, we outline convergence for a sequence of bilinear 
forms $S_n$ on certain $L^2(E,\bnu_n)$-spaces as $n\to \infty$. 
The idea comes from Mosco convergence of Dirichlet forms (U. Mosco 
\cite{Mo94}, Section 2, W. Sun, \cite{Su98}). However, we introduce 
our convergence of bilinear forms in an independent fashion. Neither 
it is formulated in the language of Dirichlet form theory, nor 
Dirichlet form theory prerequisites are required. 

Throughout the whole paper, the set of all measurable functions on a 
measurable space $S$ will be denoted by $B(S)$. If there is a notion of 
continuity, the set of all continuous functions on $S$ will be denoted 
by $C(S)$. A subscript $b$ will indicate the restriction to bounded 
functions. Similarly, we denote by $L^2$ the set of all quadratically 
integrable functions. Here we will add the space and the reference 
measure to the notation. 

\subsection{Two Classes of Bilinear Forms}

In order to introduce the basic setting, let ${\bnu}$ be a 
probability measure on a measurable space $(E,{\cal B})$ and let $( 
T_t)_{t\ge 0}$ be a strongly continuous semigroup of linear operators 
on $L^2(E,{\bnu})$. Suppose that $(T_t)_{t\ge 0}$ is associated 
with a transition probability function $P(t,x,B)$, $t\ge 0$, $x\in 
E$, $B\in {\cal B}$, i. e., $T_tf=\int f(y)\, P(t,\cdot ,dy)$, $t\ge 
0$, $f\in L^2(E,{\bnu})$. 

If we, furthermore, assume that ${\bnu}$ is an invariant measure of 
$(T_t)_{t\ge 0}$ then this means that $\int T_t f\, d{\bnu}=\int f 
\, d{\bnu}$, $t\ge 0$, $f\in L^\infty(E,{\bnu})$. Let us recall that 
the existence of an invariant probability measure ${\bnu}$ for the 
semigroup $(T_t)_{t\ge 0}$ and an associated transition probability 
function $P$ guarantee $P(t,\cdot ,E)=1$ ${\bnu}$-a.e., $t\ge 0$, and 
contractivity of $(T_t)_{t\ge 0}$ on $L^\infty (E,\bnu)$ and $L^2(E, 
{\bnu})$, cf. Lemma \ref{Lemma2.1} (c) below. 

If we do not assume that ${\bnu}$ is an invariant measure of $(T_t)_{ 
t\ge 0}$ then we suppose that $(T_t)_{t\ge 0}$ is contractive on $L^2 
(E,{\bnu})$. 
\medskip

Denoting by $(A,D(A))$ the generator of $(T_t)_{t\ge 0}$ and by 
$\langle \cdot \, , \, \cdot \rangle$ the inner product in $L^2(E, 
{\bnu})$, we introduce now the class of bilinear forms $S$ we are 
interested in. Define 
\begin{eqnarray*}
\ D(S):=\left\{u\in L^2(E,{\bnu}):\, \lim_{t\to 0}\left\langle 
\textstyle{ \frac1t}(u-T_tu)\, , \, v\right\rangle \ \mbox{\rm exists  
for all}\ v\in L^2(E,{\bnu})\right\} 
\end{eqnarray*}
and 
\begin{eqnarray*}
S(u,v):=\lim_{t\to 0}\left\langle\textstyle{ \frac1t}(u-T_tu)\, , 
\, v\right\rangle \, , \quad u\in D(S), \ v \in L^2(E,{\bnu}). 
\end{eqnarray*}
We have $D(A)=D(S)$, cf. \cite{P83}, Section 2.1, and 
\begin{eqnarray*}
S(u,v)=-\langle Au \, , \, v \rangle\, , \quad u\in D(A), \ v\in 
L^2(E,{\bnu}).
\end{eqnarray*}
Set $S(u,v):=\infty$ if $u\in L^2(E,{\bnu})\setminus D(S)$ and 
$v\in L^2(E,{\bnu})$. Let $(G_\beta)_{\beta >0}$ be the resolvent 
associated with $S$, i.e., $G_\beta =(\beta -A)^{-1}$, $\beta >0$. 

\begin{lemma}\label{Lemma2.1} 
(a) For all $u\in D(S)$, it holds that $S(u,u)\ge 0$. \\ 
(b) For all $f\in L^2(E,{\bnu})$ and all $\beta >0$, we have 
$\langle f\, , \, G_\beta f\rangle \ge 0$ where $\langle f\, , \, 
G_\beta f\rangle =0$ holds if and only if $f=0$. In addition, it 
holds that $\langle f - \beta G_\beta f\, , \, f\rangle \ge 0$. 
\\ 
(c) For all $v\in L^2(E,{\bnu})$, the function $(0,\infty )\ni 
t\to \langle T_tv\, , \, T_tv\rangle$ is decreasing. In particular, 
$(T_t)_{t\ge 0}$ is contractive on $L^2(E,{\bnu})$. 
\end{lemma} 
Proof. (a) Let us first assume that ${\bnu}$ is an invariant measure 
of $(T_t)_{t\ge 0}$. 

Let $v\in L^2(E,{\bnu})$ and $v_n:=(v\wedge n)\vee (-n)$, 
$n\in {\Bbb N}$. By the Schwarz inequality, we have $(T_tv_n)^2= 
\left(\int v_n(y)\, P(t,\cdot ,dy)\right)^2\le \int v^2_n(y)\, P(t, 
\cdot ,dy)=T_t v_n^2$, $t>0$, $n\in {\Bbb N}$. From this and the 
fact that ${\bnu}$ is an invariant measure of $(T_t)_{t\ge 0}$, 
it can be concluded that $\int (T_tv_n)^2\, d{\bnu}\le\int T_t 
v_n^2\, d{\bnu}=\int v_n^2\, {\bnu}$, $t>0$, $n\in {\Bbb N}$. 
Letting $n\to\infty$, we finally get 
\begin{eqnarray}\label{2.1}
\langle T_tv\, , \, T_tv\rangle \le \langle v\, , \, v\rangle \, , 
\quad t>0. \vphantom{\int}
\end{eqnarray}
With the Schwarz inequality, $\langle T_tv\, , \, v\rangle^2 \le 
\langle T_tv\, , \, T_tv\rangle \cdot \langle v\, , \, v\rangle$, 
this results in 
\begin{eqnarray}\label{2.2}
|\langle T_tv\, , \, v\rangle |\le \langle v\, , \, v\rangle \, , 
\quad t>0. \vphantom{\int}
\end{eqnarray}
If we do not suppose that ${\bnu}$ is an invariant measure of $(T_t 
)_{t\ge 0}$ then we get (\ref{2.2}) directly from contractivity of 
$(T_t)_{t\ge 0}$ on $L^2(E,\bnu)$. Now, $S(u,u)\ge 0$, $u\in D(S)$, 
is a consequence of the definition of $(S,D(S))$ and relation 
(\ref{2.2}). 
\medskip

\nid
(b) Let $f\in L^2(E,{\bnu})$ and $\beta >0$. Set $v:=G_\beta f$. We 
have 
\begin{eqnarray}\label{2.3}
\langle f\, , \, G_\beta f\rangle &=&\langle \beta G_\beta f - A 
G_\beta f\, , \, G_\beta f\rangle \nonumber \\ 
&=&\beta\langle v\, , \, v\rangle - \langle Av\, , \, v\rangle 
\nonumber \\ 
&=&\beta\langle v\, , \, v\rangle + S(v,v)\, . 
\end{eqnarray}
Now, the result of (a) yields $\langle f\, , \, G_\beta f\rangle 
\ge 0$. According to (a) and (\ref{2.3}), $\langle f\, , \, 
G_\beta f \rangle =0$ implies $v=0$ and thus $f=0$. Finally,  
$\langle f - \beta G_\beta f\, , \, f\rangle \ge 0$ is an 
immediate consequence of $(\ref{2.2})$. 
\medskip

\nid 
(c) This follows from relation (\ref{2.1}) \qed 
\medskip 

For later use let us introduce a second notion of a bilinear form 
for a more specified class of semigroups $(T_t)_{t\ge 0}$, satisfying 
the hypotheses of the present subsection. If we assume for a moment 
that, for every $u\in L^2(E,{\bnu})$, $t\to\left\langle\textstyle 
T_tu\, ,\, u\right\rangle$ is convex then $\lim_{t\to 0}\left 
\langle\textstyle{\frac1t}(u-T_tu)\, , \, u\right\rangle$ either 
exists or tends to $+\infty$. Let us define $D(S^c):=\left\{u\in 
L^2(E,{\bnu}):\vphantom{\textstyle{\frac1t}}\right.$ $\left.\lim_{t 
\to 0}\left\langle\textstyle{\frac1t}(u-T_tu)\, ,\, u\right\rangle< 
\infty\right\}$. From (\ref{2.2}) it can be concluded that the set 
$D(S^c)$ is linear. Furthermore, $D(S)\subseteq D(S^c)$. As in the 
parallelogram identity we obtain the existence of $\lim_{t\to 0} 
\textstyle{\frac1t}\left(\left\langle(u-T_tu)\, ,\, v\right\rangle 
+\left\langle(v-T_tv)\, ,\, u\right\rangle\right)$, $u,v\in D(S^c)$. 
Motivated by that, let us furthermore assume that $S^c(u,v):=\lim_{ 
t\to 0}\left\langle\textstyle{\frac1t}(u-T_tu)\, , \, v\right 
\rangle$ exists for all $u,v\in D(S^c)$. For $u,v\in L^2(E,{\bnu})$ 
set $S^c(u,v):=\infty$ if $u\not\in D(S^c)$ or $v\not\in D(S^c)$. 
In particular, such a form $(S^c,D(S^c))$ exists if $A$ is self-adjoint. 

\subsection{Analysis on a Sequence of $L^2$-Spaces}

From now on, suppose we are given mutually orthogonal 
probability measures ${\bnu}_n$, $n\in {\Bbb N}$, and ${\bnu}$ 
on $(E,{\cal B})$. Furthermore, suppose that ${\bnu}$ is a 
measure with countable base on $(E,{\cal B})$. In particular, 
assume that there are mutually exclusive subsets $E_n$, $n\in 
{\Bbb N}$, of $E$, such that ${\bnu}_n(E\setminus E_n)=0$. Let 
$\alpha_n$, $n\in\{0\}\cup{\Bbb N}$, 
be a sequence of positive numbers with $\sum_{n=0}^\infty\alpha_n 
=1$. Define ${\Bbb M}:=\alpha_0{\bnu}+\sum_{n=1}^\infty\alpha_n{ 
\bnu}_n$. We say that $u\in\bigcap_{n\in{\Bbb N}}L^2(E,{\bnu}_n 
)\cap L^2(E,{\bnu})$ if $u$ is an equivalence class consisting 
of all everywhere defined ${\cal B}$-measurable functions 
satisfying $f_1=f_2$ ${\Bbb M}$-a.e. if $f_1,f_2\in u$ and 
$\int u^2\, d{\bnu}_n <\infty$, $n\in {\Bbb N}$, as well as 
$\int u^2\, d{\bnu}<\infty$. Let $\langle\cdot \, , \, \cdot 
\rangle_n$ denote the inner product in $L^2(E,{\bnu}_n)$, $n\in 
{\Bbb N}$, and let $\langle \cdot \, , \, \cdot \rangle$ denote 
the inner product in $L^2(E,{\bnu})$. Introduce 
\begin{eqnarray*}
{\cal D}:=\left\{u\in \bigcap_{\, n\in {\Bbb N}}L^2(E,{\bnu}_n) 
\cap L^2 (E,{\bnu}): \langle u\, , \, u\rangle_n\stack{n\to\infty} 
{\lra}\langle u\, , \, u\rangle \right\}
\end{eqnarray*}
and suppose that the set ${\cal D}$ contains a linear subset 
${\cal C}$. Also introduce  
\begin{itemize}
\item[(${\cal C}1$)] ${\cal C}$ is dense in $L^2(E,{\bnu})$. 
\end{itemize}
\begin{definition}\label{Definition2.2} 
(a) A sequence $\vp_n\in {\cal C}$, $n\in {\Bbb N}$, is said to be 
{\it w-convergent} to $\vp\in L^2(E,{\bnu})$ as $n\to\infty$ if 
\begin{itemize} 
\item[(i)] $\langle\vp_n\, , \, \psi\rangle_n\stack{n\to 
\infty}{\lra}\langle\vp\, , \, \psi\rangle$ for all $\psi\in 
{\cal C}$.  
\end{itemize} 
(b) A sequence $\psi_n\in {\cal C}$, $n\in {\Bbb N}$, is said to 
be {\it s-convergent} to $\psi\in L^2(E,{\bnu})$ as $n\to 
\infty$ if 
\begin{itemize} 
\item[(i)] $\psi_n$ $w$-converges to $\psi$ as $n\to\infty$ and 
\item[(ii)] $\langle \psi_n\, , \, \psi_n\rangle_n \stack{n\to 
\infty}{\lra}\langle \psi\, , \, \psi\rangle$.
\end{itemize} 
(c) $w$-convergence or $s$-convergence of subsequences 
$\vp_{n_k}\in {\cal C}$ or $\psi_{n_k}\in {\cal C}$, respectively, 
will mean that in (a) or (b) the index $n\in {\Bbb N}$ is replaced 
with $n_k\in {\Bbb N}$. 
\end{definition}
{\bf Remark} (1) Note that according to the definition of ${\cal 
C}$, we have the following implication: If $\psi_n\in {\cal C}$, 
$n\in {\Bbb N}$, $s$-converges to $\psi\in {\cal C}$ then 
\begin{eqnarray*}
\langle\psi_n-\psi \, , \, \psi_n-\psi\rangle_n &=&\langle\psi_n 
\, , \, \psi_n\rangle_n -2\langle\psi_n\, , \, \psi\rangle_n + 
\langle\psi \, , \, \psi\rangle_n \\ 
&\stack{n\to\infty}{\lra}& 0\, . 
\end{eqnarray*}
Also we observe that for $\psi\in {\cal C}$ the sequence $\psi_n 
:=\psi$, $n\in {\Bbb N}$, $s$-converges to $\psi$ as $n\to\infty$. 
\begin{proposition}\label{Proposition2.3} 
Suppose that condition (${\cal C}1$) is satisfied. 
(a) Let $\vp_n\in {\cal C}$, $n\in {\Bbb N}$, be a sequence such 
that $\langle \vp_n\, , \, \vp_n\rangle_n$ is bounded. Then there 
exists a subsequence $\vp_{n_k}\in {\cal C}$, $k\in {\Bbb N}$, 
$w$-convergent to some $\vp\in L^2(E,{\bnu})$ as $k\to\infty$. \\  
(b) Let $\vp_n\in {\cal C}$, $n\in {\Bbb N}$, be a sequence 
$w$-convergent to $\vp\in L^2(E,{\bnu})$. Then we have 
\begin{eqnarray*}
\liminf_{n\to\infty}\langle\vp_n \, ,\, \vp_n\rangle_n\ge \langle 
\vp\, , \, \vp\rangle\, . 
\end{eqnarray*}
\end{proposition}
Proof. {\it Step 1 } The space $L^2(E,{\bnu})$ is separable. 
Thus, there is a sequence $w_n\in L^2(E,{\bnu})$, $n\in {\Bbb 
N}$, such that for every $w\in L^2(E,{\bnu})$, there is a 
subsequence $w_{n_k}\in L^2(E,{\bnu})$, $k\in {\Bbb N}$, with 
$w_{n_k}\stack{k\to \infty}{\lra}w$ in $L^2(E,{\bnu})$. Since 
${\cal C}$ is dense in $L^2(E,{\bnu})$ (condition (${\cal C}1$)) 
there exists a sequence $v_s^t\in {\cal C}$, $s,t\in {\Bbb 
N}$, such that $\left\langle w_s-v_s^t\, , \, w_s-v_s^t\right 
\rangle_n<\frac1t$. Let $v_n\in {\cal C}$, $n\in {\Bbb N}$, be a 
sequence with $\{v_n:n\in {\Bbb N}\}=\{v_s^t:s,t\in {\Bbb N}\}$. 
Then $v_n\in {\cal C}$, $n\in {\Bbb N}$, is a sequence such that 
for every $v\in {\cal C}$, there is a subsequence $v_{n_k}\in 
{\cal C}$, $k\in {\Bbb N}$, with $v_{n_k}\stack{k\to\infty} 
{\lra}v$ in the norm of $L^2(E,{\bnu})$. Consequently, for 
$\ve>0$ and $\delta >0$, there is an $r\in {\Bbb N}$, such 
that 
\begin{eqnarray}\label{2.4}
\langle v_r-v\, , \, v_r-v\rangle^{1/2} <\delta\cdot\ve  
\end{eqnarray}
(in Step 3 we will use the $\delta$ in this notation). Recalling 
the definition of ${\cal C}$, we verify now the existence of $n_0 
\in {\Bbb N}$ such that 
\begin{eqnarray}\label{2.5}
\langle v_r-v\, , \, v_r-v\rangle^{1/2}_n <\ve\, , \quad n>n_0.  
\end{eqnarray}
{\it Step 2 } Let $v_k\in {\cal C}$, $k\in {\Bbb N}$, be the 
sequence introduced in Step 1 and let ${\cal L}$ denote the set 
of all finite linear combinations of ${v}_k$, $k\in {\Bbb N}$. 
Without loss of generality, we may assume that $|\langle\vp_n\, 
, \, \vp_n\rangle_n|\le 1$, $n\in {\Bbb N}$. Note that ${\cal L} 
\subseteq{\cal C}$ and that, for all ${v}\in {\cal L}$, we have 
\begin{eqnarray}\label{2.6}
|\langle\vp_n\, , \, {v}\rangle_n|\le \langle{v}\, , \, {v}
\rangle^{1/2}_n\, , \quad n\in {\Bbb N}. 
\end{eqnarray}
From $\vp_n$, $n\in {\Bbb N}$, extract a subsequence $\vp_{ 
n_{k_1}}$, $k_1\in {\Bbb N}$, such that $\langle\vp_{n_{k_1}} 
\, , \, {v}_1\rangle_{n_{k_1}}$ converges to some $\gamma 
({v}_1)$ as $k_1\to\infty$ and proceed as follows: From 
$\vp_{n_{k_m}}$, $k_m\in {\Bbb N}$, extract a subsequence 
$\vp_{n_{k_{m+1}}}$, $k_{m+1}\in {\Bbb N}$, such that $\langle 
\vp_{n_{k_{m+1}}} \, , \, {v}_{m+1}\rangle_{n_{k_{m+1}}}$ 
converges to some $\gamma ({v}_{m+1})$ as $k_{m+1}\to\infty$, 
$m\in {\Bbb N}$. For an arbitrary ${v}\in {\cal L}$, say ${v} 
:=\sum_{i=1}^ra_i{v}_{k_i}$, $a_1,\ldots ,a_r\in {\Bbb R}$, 
set $\gamma({v}):=\sum_{i=1}^ra_i\gamma({v}_{k_i})$. Let 
$l_m$ be the $m$-th member of $n_{k_m}$, $k_m\in {\Bbb N}$. 
According to the above selection procedure, we obtain 
\begin{eqnarray}\label{2.7}
\langle\vp_{l_m}\, , \, {v}\rangle_{l_m}\stack{m\to\infty} 
{\lra}\gamma ({v})\, , \quad {v}\in {\cal L}. 
\end{eqnarray}
Now, (\ref{2.6}), (\ref{2.7}), and $v\in {\cal L}\subseteq 
{\cal C}\subseteq {\cal D}$ imply $|\gamma ({v})|\le \langle{v} 
\, , \, {v}\rangle^{1/2}$, $v\in {\cal L}$, 
which means that $\gamma$ defines a bounded linear functional 
on ${\cal L}$. By (\ref{2.4}), $\{v_n:n\in {\Bbb N}\}$ is 
dense in ${\cal C}$ w.r.t. the norm in $L^2(E,{\bnu})$ and 
by condition (${\cal C}1$), ${\cal L}$ is thus dense in 
$L^2(E,{\bnu})$. Consequently, $\gamma$ can continuously be 
extended to a bounded linear functional on $L^2(E,{\bnu})$.   
Let this extension also be denoted by $\gamma$. 
There exists a $\vp\in L^2(E,{\bnu})$ such that $\gamma 
(v)=\langle\vp\, , \, v\rangle$, $v\in L^2(E,{\bnu})$, and 
with (\ref{2.7}), we have 
\begin{eqnarray}\label{2.8}
\langle\vp_{l_m}\, , \, {v}\rangle_{l_m}\stack{m\to\infty} 
{\lra}\langle\vp\, , \, {v}\rangle \, , \quad {v}\in {\cal L}. 
\end{eqnarray}
\nid
{\it Step 3 } Our goal is now to demonstrate that (\ref{2.8}) 
holds true for all $v\in{\cal C}$. To this end, let $\ve >0$ 
and let $v\in {\cal C}$. Again, let $v_k\in {\cal C}$, $k\in 
{\Bbb N}$, be the sequence introduced in Step 1. Then there 
exist $r\in {\Bbb N}$ and $n_0\in {\Bbb N}$ with $\langle{v} 
-{v}_r\, , \, {v}-{v}_r\rangle^{1/2}_n <\ve$ for all $n>n_0$  
and $\langle{v}-{v}_r\, , \, {v} -{v}_r\rangle^{1/2}<\ve/\langle 
\vp\, , \, \vp\rangle^{1/2}$, cf. (\ref{2.4}) and (\ref{2.5}). 
In addition, there is an $N\in {\Bbb N}$ such that $|\langle 
\vp_{l_m}\, , \, {v}_r\rangle_{l_m}-\langle\vp\, , \, {v}_r 
\rangle|<\ve$ for all $m>N$, recall (\ref{2.8}). For $m>N$ and 
$l_m>n_0$, we thus have 
\begin{eqnarray*}
&&|\langle\vp_{l_m}\, , \, {v}\rangle_{l_m}-\langle\vp\, , \, 
{v}\rangle| \\ 
&&\hspace{1.5cm}\le |\langle\vp_{l_m}\, , \, {v}\rangle_{l_m}- 
\langle\vp_{l_m} \, , \, {v}_r\rangle_{l_m}| + |\langle\vp_{ 
l_m}\, , \, {v}_r\rangle_{l_m}-\langle\vp\, , \, {v}_r\rangle| 
+|\langle\vp\, , \, {v}_r\rangle-\langle\vp\, , \, {v}\rangle| 
\\
&&\hspace{1.5cm}\le\langle\vp_{l_m}\, , \, \vp_{l_m}\rangle^{1/2 
}_{l_m}\langle {v} -{v}_r\, , \, {v} -{v}_r\rangle^{1/2}_{ 
l_m}+\ve + \langle\vp\, , \, \vp\rangle^{1/2}\langle {v}-{v}_r
\, , \, {v} -{v}_r\rangle^{1/2} \\ 
&&\hspace{1.5cm}<3\ve\, , 
\end{eqnarray*}
where the last inequality uses the assumption $|\langle\vp_n 
\, , \, \vp_n\rangle_n|\le 1$, $n\in {\Bbb N}$. Therefore, 
\begin{eqnarray*}
\langle\vp_{l_m}\, , \, {v}\rangle_{l_m}\stack{m\to\infty}{\lra}
\langle\vp\, , \, {v}\rangle\, , \quad {v}\in  {\cal C}\, . 
\end{eqnarray*}
Hence, we have shown that the subsequence $\vp_{l_m}\in {\cal C}$, 
$m\in {\Bbb N}$, of $\vp_n\in {\cal C} $, $n\in {\Bbb N}$, is 
$w$-convergent to $\vp \in L^2(E,{\bnu})$. 
\medskip

\nid 
(b) Let $\ve >0$ and $\t \vp\in {\cal C}$ such that $\langle\t \vp 
-\vp \, , \, \t \vp-\vp\rangle^{1/2} <\ve$, cf. condition (${\cal 
C}1$). Since $\vp_n$, $n\in {\Bbb N}$, $w$-converges to $\vp\in 
{\cal C}$, it follows that  
\begin{eqnarray*}
\langle\vp_n\, , \, \vp_n\rangle_n-\langle\vp_n-\t \vp \, , \, 
\vp_n-\t \vp\rangle_n&=&2\langle\vp_n\, , \, \t \vp\rangle_n
-\langle \t \vp\, , \, \t \vp\rangle_n \\ 
&\stack{n\to\infty}{\lra}&2\langle \vp\, , \, \t \vp\rangle
-\langle \t \vp\, , \, \t \vp\rangle \\ 
&=&\langle\vp\, , \, \vp\rangle + 2\langle\vp\, , \, \t \vp-\vp 
\rangle + \left(\langle\vp\, , \, \vp\rangle -\langle\t \vp\, , \, 
\t \vp\rangle\right)\, .  
\end{eqnarray*}
From $|\langle\vp\, , \, \t \vp-\vp\rangle|\le\langle\vp\, , \, 
\vp\rangle^{1/2}\cdot \ve$ and $|\langle\vp\, , \, \vp\rangle 
-\langle\t \vp\, , \, \t \vp\rangle|=|\langle\vp +\t \vp\, , \, 
\vp -\t \vp\rangle|\le 2\langle\vp\, , \, \vp\rangle^{1/2}\cdot 
\ve +\ve^2$ we finally get $\liminf_{n\to\infty}\langle\vp_n\, , 
\, \vp_n\rangle_n\ge\langle \vp\, , \, \vp\rangle$. 
\qed 

\subsection{Convergence of Bilinear Forms}

Throughout the paper, for every $n\in {\Bbb N}$, let $S_n$ be a 
bilinear form on $L^2(E,{\bnu}_n)$, and let $S$ be a bilinear form 
on $L^2(E,{\bnu})$. Suppose that ${\bnu}_n$ is an invariant measure 
of the strongly continuous semigroup $(T_{n,t})_{t\ge 0}$ in $L^2(E, 
{\bnu}_n)$ associated with $S_n$, $n\in {\Bbb N}$. Furthermore, 
suppose that $(T_{n,t})_{t\ge 0}$ possesses a transition probability 
function, $n\in {\Bbb N}$. Also assume that we are given a semigroup 
$(T_t^b)_{t\ge 0}$ in $B_b(E)$ possessing a transition probability 
function by means of which $(T_t^b)_{t\ge 0}$ induces a strongly 
continuous contraction semigroup $(T_t)_{t\ge 0}$ in $L^2(E,{\bnu})$. 
If no ambiguity is possible we will drop the superscript $b$ from the 
notation. Suppose that $(T_t)_{t\ge 0}$ is associated with $S$. Note  
that ${\bnu}$ is not necessarily an invariant probability measure of 
the semigroup $(T_t)_{t\ge 0}$. Note furthermore, that Lemma 
\ref{Lemma2.1} (a) guarantees non-negativity of the bilinear forms 
$S_n$, $n\in {\Bbb N}$, and $S$. 

Furthermore, let $G_{n,\beta}$ and $G_\beta$, $\beta >0$, denote the 
families of resolvents associated with $(T_{n,t})_{t\ge 0}$ and $(T_{t} 
)_{t\ge 0}$ and let $A_n$ and $A$ be the generators of the semigroups 
$(T_{n,t})_{t\ge 0}$, $n\in {\Bbb N}$, and $(T_{t})_{t\ge 0}$. 
Introduce  
\begin{itemize}
\item[(${\cal C}2$)] ${\cal G}:=\{G_\beta g:g\in {\cal C}\cap B_b(E), 
\ \beta >0\}\subseteq {\cal C}$ and ${\cal G}_n:=\{G_{n, \beta}g:g 
\in {\cal C}\cap B_b(E),\ \beta >0\}\subseteq {\cal C}$, $n\in {\Bbb 
N}$, in the sense that for every $g\in {\cal C}\cap B_b(E)$ and 
$\beta>0$, there is a $u\in {\cal C}$ with $G_\beta g=u$ ${\bnu} 
$-a.e. and furthermore, for every $n\in {\Bbb N}$, there exists a $v_n 
\in {\cal C}$ such that $G_{n,\beta}g=v_n$ ${\bnu}_n$-a.e.  
\end{itemize}
For $\beta >0$, $n\in {\Bbb N}$, $\vp_n\in D(S_n)$, $\psi_n\in 
L^2(E,{\bnu}_n)$, set $S_{n,\beta}(\vp_n,\psi_n):=\beta\langle \vp_n 
\, , \, \psi_n \rangle_n + S_n(\vp_n,\psi_n)$, and for $\vp\in D(S)$, 
$\psi\in L^2(E,{\bnu})$, define $S_{\beta}(\vp,\psi):=\beta\langle 
\vp\, , \, \psi\rangle +S(\vp,\psi)$. Furthermore, set $S_{n,\beta} 
(\vp_n,\psi_n)=\infty$ if $S_n(\vp_n,\psi_n)=\infty$ and $S_{\beta} 
(\vp,\psi)=\infty$ if $S(\vp,\psi)=\infty$. 
\begin{definition}\label{Definition2.4} 
(a) We say that $S_n$, $n\in {\Bbb N}$, {\it pre-converges} to $S$ if 
we have the following.
\begin{itemize} 
\item[(i)] For every $\vp\in L^2(E,{\bnu})$, every subsequence $n_k$, 
$k\in {\Bbb N}$, of indices, and every subsequence $\vp_{n_k}\in {\cal 
C}$, $k\in {\Bbb N}$, $w$-converging to $\vp$, we have 
\begin{eqnarray*}
S(\vp,\vp )\le \liminf_{k\to\infty}S_{n_k}(\vp_{n_k},\vp_{n_k})\, . 
\end{eqnarray*}
\item[(ii)] For every $\psi\in D(S)$, there exists a sequence $\psi_n 
\in D(S_n)\cap {\cal C}$, $n\in {\Bbb N}$, $s$-converging to $\psi$ 
such that 
\begin{eqnarray*}
\limsup_{n\to\infty} S_n(\psi_n ,\psi_n )\le S(\psi ,\psi )\, .  
\end{eqnarray*}
\end{itemize} 
\end{definition}
\medskip 

\nid 
{\bf Remarks} (2) In symmetric Dirichlet form theory, conditions 
(i) and (ii) are known as {\it Mosco convergence}, cf. \cite{Mo94}. 
\medskip 

\nid 
(3) Imposing condition (i) on $S_n$, $n\in {\Bbb N}$, and $S$, and 
assuming $({\cal C}2)$, we 
implicitly require that, for all $\beta >0$ and $g\in {\cal C}\cap 
B_b(E)$, every $w$-limit $\t u$ of a sequence $u_n:=G_{n,\beta}g$, 
$n\in {\Bbb N}$, belongs to $D(S)$: Recalling the definition of $(S, 
D(S))$ and condition $({\cal C}2)$, this can be verified by $S(\t u, 
\t u)+\beta\langle\t u\, , \, \t u \rangle \le \liminf_{n\to \infty} 
(S_n(u_n,u_n) +\beta\langle u_n\, , \, u_n\rangle_n)=\liminf_{n\to 
\infty}\langle g\, , \, u_n\rangle_n\le\frac{1}{\beta}\sup_{x\in E} 
g(x)^2 <\infty$.  

\begin{lemma}\label{Lemma2.5} 
Let $S_n$, $n\in {\Bbb N}$, be pre-convergent to $S$. Furthermore,  
let $w_n\in {\cal C}$, $n\in {\Bbb N}$, be a sequence $w$-converging 
to  some $w\in D(S)$. Finally, let $v_n\in D(S_n)\cap {\cal C}$, 
$n\in {\Bbb N}$, be a sequence $s$-converging to $v\in D(S)$ 
satisfying condition (ii) of Definition \ref{Definition2.4}. \\ 
(a) If $\, \limsup_{n\to\infty}S_n(w_n,w_n)<\infty$ then the limit 
$\lim_{n\to\infty}(S_n(v_n,w_n)+S_n(w_n,v_n))$ exists and we have 
\begin{eqnarray*}
\lim_{n\to\infty}(S_n(v_n,w_n)+S_n(w_n,v_n))=S(v,w)+S(w,v)\, . 
\end{eqnarray*}
(b) If $\beta >0$ and $\, \limsup_{n\to\infty}S_{n,\beta}(w_n,w_n)< 
\infty$ then 
\begin{eqnarray*}
\lim_{n\to\infty}(S_{n,\beta}(v_n,w_n)+S_{n,\beta}(w_n,v_n))= 
S_\beta (v,w)+S_\beta (w,v)\, . 
\end{eqnarray*}
The lemma holds also for subsequences $n_k$, $k\in {\Bbb N}$, of 
indices. 
\end{lemma}
Proof. Under the above assumptions, we have 
\begin{eqnarray*}
&&\hspace{-1.5cm}-\limsup_{n\to\infty}S_n(w_n,w_n)+S(v+w,v+w)-S(v,v)\\ 
&\le&\liminf_{n\to\infty}\left\{-S_n(w_n,w_n)+S_n(v_n+w_n,v_n+w_n)- 
S_n(v_n,v_n)\right\} \\  
&=&\liminf_{n\to\infty}\left\{S_n(v_n,w_n)+S_n(w_n,v_n)\right\} \\  
\end{eqnarray*}
as well as 
\begin{eqnarray*}
&&\hspace{-1.5cm}\limsup_{n\to\infty}\left\{S_n(v_n,w_n)+S_n(w_n,v_n) 
\right\} \\  
&=&\limsup_{n\to\infty}\left\{S_n(w_n,w_n)-S_n(v_n-w_n,v_n-w_n)+ 
S_n(v_n,v_n)\right\} \\  
&\le&\limsup_{n\to\infty}S_n(w_n,w_n)-S(v-w,v-w)+S(v,v)\, . 
\end{eqnarray*}
Thus, 
\begin{eqnarray*}
\limsup_{n\to\infty}\left|S_n(v_n,w_n)+S_n(w_n,v_n)-S(v,w)-S(w,v) 
\right|\le\limsup_{n\to\infty}S_n(w_n,w_n)-S(w,w)\, .  
\end{eqnarray*}
For $\ve >0$, this relation implies 
\begin{eqnarray*}
&&\limsup_{n\to\infty}\left|S_n(v_n,w_n)+S_n(w_n,v_n)-S(v,w)-S(w,v) 
\right|\\ 
&&\hspace{0.8cm}=\limsup_{n\to\infty}\left|S_n(v_n/\ve,\ve w_n)+S_n 
(\ve w_n,v_n/\ve)-S(v/\ve,\ve w)-S(\ve w,v/\ve)\right| \\ 
&&\hspace{0.8cm}\le\limsup_{n\to\infty}S_n(\ve w_n,\ve w_n)-S(\ve w, 
\ve w) \\ 
&&\hspace{0.8cm}=\ve^2\left\{\limsup_{n\to\infty}S_n(w_n,w_n)-S(w,w) 
\right\} \, .  
\end{eqnarray*}
Recalling that $\ve >0$ has been chosen arbitrarily and that, 
according to Definition \ref{Definition2.4} and the above assumptions 
of this lemma, we have $S(w,w)\le \limsup_{n\to\infty}S_n(w_n,w_n) 
<\infty$, we obtain $\lim_{n\to\infty}(S_n(v_n,w_n)+S_n(w_n,v_n)) 
=S(v,w)+S(w,v)$. The second statement follows from similar arguments 
since (ii) of Definition \ref{Definition2.2} (b) and Proposition 
\ref{Proposition2.3} (b) imply that conditions (i) and (ii) of 
Definition \ref{Definition2.4} also hold for $S_{n,\beta}$ and 
$S_\beta$ instead of $S_n$, $n\in {\Bbb N}$, and $S$. 
\qed
\begin{proposition}\label{Proposition2.6} 
Assume $({\cal C}1)$ and $({\cal C}2)$. 
Let $S_n$, $n\in {\Bbb N}$, be pre-convergent to $S$. Furthermore, 
let $\beta >0$, $g\in {\cal C}\cap B_b(E)$, $u_n:=G_{n,\beta}g$, 
$n\in {\Bbb N}$, and $u:=G_{\beta}g$. \\ 
(a) For every subsequence $u_{n_k}$, $k\in {\Bbb N}$, there exists 
another subsequence $w$-converging to some $\t u\in D(S)$. \\ 
(b) The following are equivalent. 
\begin{itemize}
\item[(iii)] For every $\psi\in D(S)$ and every sequence $\psi_n\in 
D(S_n)\cap {\cal C}$, $n\in {\Bbb N}$, $s$-converging to $\psi$ such 
that (ii) of Definition \ref{Definition2.4} is satisfied, we have 
\begin{eqnarray*}
\lim_{n\to\infty}S_{n,\beta}(\psi_n,G_{n,\beta}g)=S_\beta(\psi, 
G_\beta g)\, . 
\end{eqnarray*}
\item[(iii')] For all $\t u\in D(S)$ such that there is a 
subsequence $u_{n_k}$, $k\in {\Bbb N}$, $w$-converging to $\t u$, 
we have $S_\beta (\psi, u)+S_\beta (u,\psi)=S_\beta (\psi ,\t u) 
+S_\beta (\t u,\psi )$, $\psi\in D(S)$. 
 
\item[(iv)] For all $\t u\in D(S)$ such that there is a 
subsequence $u_{n_k}$, $k\in {\Bbb N}$, $w$-converging to $\t u$, 
we have $\t u=u$ ${\bnu}$-a.e. 
\end{itemize} 
\end{proposition}
Proof. {\it Step 1 } Part (a) follows from $\langle u_n\, , \, u_n 
\rangle_n^{1/2}\le\frac1\beta \langle g\, , \, g\rangle_n^{1/2}\le 
\frac1\beta\sup_{x\in E}|g(x)|$, $n\in {\Bbb N}$, condition $({\cal 
C}2)$, Proposition \ref{Proposition2.3} (a), and Remark (3). 
\medskip 

\nid
{\it Step 2 } (iii') implies (iii): Recall Lemma 
\ref{Lemma2.5} and note that $\limsup_{n\to\infty}S_{n,\beta}(u_n, 
u_n)=\limsup_{n\to\infty}\langle g\, , \,G_{n,\beta}g\rangle_n\le 
\frac{1}{\beta}\sup_{x\in E}g(x)^2<\infty$. Turning to subsequences 
if necessary and keeping part (a) in mind, we may state 
\begin{eqnarray*}
\lim_{n\to\infty}S_{n,\beta}(\psi_n,G_{n,\beta}g) &=&\lim_{n\to\infty}
\left(S_{n,\beta}(\psi_n,G_{n,\beta}g)+S_{n,\beta}(G_{n,\beta}g, 
\psi_n)\right)-\lim_{n\to\infty}S_{n,\beta}(G_{n,\beta}g,\psi_n) \\ 
&=&\left(S_\beta (\psi,\t u)+S_\beta(\t u,\psi)\right)-\lim_{n\to 
\infty}\langle g\, , \, \psi_n\rangle_n \\ 
&=&S_\beta (\psi,\t u)+S_\beta(\t u,\psi)-\langle g\, , \, \psi\rangle 
\\ 
&=&S_\beta (\psi,\t u)+S_\beta(\t u,\psi)-S_\beta (u,\psi) \\ 
&=&S_\beta (\psi,u) \\ 
&=&S_\beta (\psi,G_\beta g)\, .   
\end{eqnarray*}
Since the right-hand side is independent of the possible choice of a 
subsequence, the limit (iii) exists. 
\medskip 

\nid
{\it Step 3 } (iii) implies (iii'): This becomes evident after 
rearranging the chain of equations in Step 2. 
\medskip 

\nid
{\it Step 4 } (iii') implies (iv): Applying (iii') to $\psi:=\t u$ 
as well as to $\psi:=u$, we obtain $S_\beta (\t u-u,\t u-u)=0$. This 
yields (iv). 
\medskip 

\nid
{\it Step 5 } That (iv) implies (iii') is trivial. 
\qed 
\medskip 

\nid
{\bf Remark} (4) In case of symmetric forms $S_n$, $n\in {\Bbb N}$, and 
$S$ (i.e., $S_n(\vp_n,\psi_n)=S_n(\psi_n,\vp_n)$, $n\in {\Bbb N}$, $\vp_n, 
\psi_n\in D(S_n)$ and $S(\vp, \psi)=S(\psi,\vp)$, $\vp,\psi\in D(S)$), 
conditions (iii), (iii'), and (iv) are trivial: This follows from 
\begin{eqnarray*}
\lim_{n\to\infty}S_{n,\beta}(\psi_n,G_{n,\beta}g)&=&\lim_{n\to\infty} 
S_{n,\beta}(G_{n,\beta}g,\psi_n) \\ 
&=&\lim_{n\to\infty}\langle g\, , \, \psi_n\rangle_n \\ 
&=&\langle g\, , \, \psi\rangle \\ 
&=&S_\beta(G_\beta g,\psi) \\ 
&=&S_\beta(\psi,G_\beta g)\, . 
\end{eqnarray*}
{\bf Definition\ \ref{Definition2.4}\ continued}  
(a) We say that $S_n$, $n\in {\Bbb N}$, {\it converges} to $S$ if 
we have (i), (ii), and (iii). 
\medskip 

\nid 
{\bf Remark} (5) Among the equivalent conditions to be added 
to (i) and (ii) in order to handle non-symmetry of the forms $S_n$, 
$n\in {\Bbb N}$, we have selected (iii) since we have used it in 
Remark (4) to verify validity of (iii), (iii'), and (iv) in the case 
of symmetry. 
\medskip 

Our objective is now to demonstrate that the above notion of 
convergence of forms $S_n$ to the form $S$ as $n\to\infty$ is 
sufficient for $s$-convergence of resolvents. 
\begin{theorem}\label{Theorem2.7} 
Suppose that conditions (${\cal C}1$) and (${\cal C}2$) are satisfied 
and that $S_n$, $n\in {\Bbb N}$, converges to $S$ in the sense of 
Definition \ref{Definition2.4}. Then, for all $g\in {\cal C}\cap 
B_b(E)$ and $\beta >0$, $G_{n,\beta}g$ $s$-converges to $G_\beta g$ 
as $n\to\infty$.  
\end{theorem} 
Proof. {\it Step 1 } Fix $g\in {\cal C}\cap B_b(E)$ and $\beta >0$. Set 
$u_n:=G_{n,\beta}g$. Because of Proposition \ref{Proposition2.6} (a), 
there exists a subsequence $u_{n_k}$, $k\in {\Bbb N}$, $w$-converging 
to some $\t u\in D(S)$. Let $u:=G_\beta g$. 

From Proposition \ref{Proposition2.6} (b) it follows that $\t u=G_\beta 
g=u$. Thus, $u_n=G_{n,\beta}g$, $n\in {\Bbb N}$, $w$-converges to $u= 
G_\beta g$, independent of the possible choice of a subsequence above. 
\medskip

\nid 
{\it Step 2 } It remains to show that $\langle u_n\, , \, u_n 
\rangle_n \stack{n\to\infty}{\lra}\langle u \, , \, u \rangle$. 
Recalling condition $({\cal C}2)$, we figure 
\begin{eqnarray*}
\beta \langle u\, , \, u\rangle + S(u,u)&=&\langle g\, , \, u\rangle \\ 
&=&\lim_{n\to\infty}\langle g\, , \, u_n\rangle_n \\ 
&=&\lim_{n\to\infty}\left\{\beta\langle u_n \, , \, u_n\rangle_n+ 
S_n(u_n,u_n)\right\}\, . 
\end{eqnarray*}
From this equality and Proposition \ref{Proposition2.3} (b) as well 
as Definition \ref{Definition2.4} (i), we finally derive 
$\lim_{n\to\infty} S_n(u_n,u_n)=S(u,u)$ and the desired relation 
$\lim_{n\to\infty}\langle u_n\, , \, u_n \rangle_n =\langle u 
\, , \, u \rangle$. \qed 
\medskip 

\nid
{\bf Remarks} (6) Let $u\in {\cal C}$ such that $G_\beta g=u$
$\bnu$-a.e., cf. condition (${\cal C}2$). By virtue of Theorem 
\ref{Theorem2.7} and Remark (1) we have 
\begin{eqnarray}\label{2.9}
\langle G_{n,\beta}g - u\, , \, G_{n,\beta}g - u\rangle_n \stack{n 
\to\infty}{\lra} 0\, . 
\end{eqnarray}
(7) Following the proofs from Lemma \ref{Lemma2.5} on it turns out 
that there is another version of Theorem \ref{Theorem2.7}. Instead 
of (ii) and (iii) let us require the following. 
\begin{itemize} 
\item[(ii')] For every subsequence of indices $n_q$, $q\in {\Bbb N}$, 
and $\psi\in D(S)$ there exists another subsequence $n_r$, $r\in {\Bbb N}$, 
of $n_q$, $q\in {\Bbb N}$, and $\psi_{n_r}\in D(S_{n_r})\cap {\cal C}$, 
$r\in {\Bbb N}$, $s$-converging to $\psi$ such that  
\begin{eqnarray*}
\limsup_{r\to\infty} S_{n_r}(\psi_{n_r} ,\psi_{n_r} )\le S(\psi ,\psi ) 
\, .  
\end{eqnarray*}
\item[(iii'')] For every $\psi\in D(S)$ and every subequence $\psi_{n_r} 
\in D(S_{n_r})\cap {\cal C}$, $r\in {\Bbb N}$, $s$-converging to $\psi$ 
such that (ii'), we have 
\begin{eqnarray*}
\lim_{r\to\infty}S_{n_r,\beta}(\psi_{n_r},G_{n_r,\beta}g)=S_\beta(\psi, 
G_\beta g)\, . 
\end{eqnarray*}
\end{itemize} 
Let us say that {\it $S_n$ converges to $S$ in the sense of Remark (7) of 
Section 2} if we have (i) of Definition \ref{Definition2.4}, (ii'), and 
(iii''). The version of Theorem \ref{Theorem2.7} we just have established 
reads as follows. 

{\it Suppose (${\cal C}1$), (${\cal C}2$), and assume this convergence 
of $S_n$ to $S$. Then, for all $g\in {\cal C}\cap B_b(E)$ and $\beta >0$, 
$G_{n,\beta}g$ $s$-converges to $G_\beta g$ as $n\to\infty$. } 
\medskip 

\nid
(8) Let us assume that for the limiting semigroup, $t\to\left\langle\textstyle 
T_tu\, ,\, u\right\rangle$ is convex for every $u\in L^2(E,{\bnu})$ and 
$\lim_{t\to 0}\left\langle\textstyle{\frac1t}(u-T_tu)\, , \, v\right\rangle$ 
exists for all $u,v\in D(S^c)$, cf. end of Subsection 2.1. Then the whole 
analysis of this subsection remains valid if $(S,D(S))$ is replaced by 
$(S^c,D(S^c))$. 
\medskip 

\nid
(9) Let $\1$ denote the function constant one on $E$. Assume $\1\in {\cal 
C}$. As a consequence of Theorem \ref{Theorem2.7}, $\langle\1\, ,\, \beta 
G_\beta g\rangle=\lim_{n\to\infty}\langle\1\, ,\, \beta G_{n,\beta} g 
\rangle_n=\lim_{n\to\infty}\langle\1\, ,\, g\rangle_n=\langle\1\, ,\, g 
\rangle$, $g\in {\cal C}\cap B_b(E)$, $\beta >0$, i. e., $\bnu$ is an 
invariant measure of the semigroup $(T_t)_{t\ge 0}$.

\section{Weak Convergence of Particle Processes}
\setcounter{equation}{0}

Let $E$ be a compact metric space. Let $(T_{n,t})_{t\ge 0}$, $n\in 
{\Bbb N}$, and $(T^b_t)_{t\ge 0}$ be associated with a cadlag stochastic 
processes. Our goal is to establish weak convergence of these processes 
in the Skorohod space $D_E[0,\infty)$ if the initial distributions are 
the invariant measures $\bnu_n$, $n\in {\Bbb N}$, and $\bnu$ or from a 
certain class of its perturbations. In order to be consistent with the 
preceding (sub)sections, we will keep on writing $C_b(E)$ for $C(E)$. 
\medskip

In this introductory part of the section, let us get the idea of how 
we use Theorem \ref{Theorem2.7} in the paper. For any accumulation point 
$\t{\bnu}$ of $({\bnu}_n)_{n\in {\Bbb N}}$ introduce the following. Let 
$(T_t)_{t\ge 0}$ be the semigroup in $L^2(E,\t {\bnu})$ induced by $( 
T_t^b)_{t\ge 0}$. Assume that $(T_t)_{t\ge 0}$ is strongly continuous 
and contractive on $L^2(E,\t {\bnu})$. Relative to the measure $\t\bnu$ 
introduce $A_{\t \sbnu},S_{\t\sbnu}, S^c_{\t \sbnu}$ along the lines of 
Subsection 2.1. 
\medskip 

Let us for this introduction assume that we are given a subsequence of 
indices $n_k$, $k\in {\Bbb N}$, such that ${\bnu}_{n_k}\stack{k\to 
\infty}{\Ra}\t\bnu$. Also suppose (${\cal C}1$), (${\cal C}2$), and that 
the hypotheses of Theorem  \ref{Theorem2.7} are satisfied, i. e., that 
$S_{n_k}$, $k\in {\Bbb N}$, converges to $S_{\t\sbnu}$ or $S^c_{\t\sbnu 
}$ in the sense of Definition \ref{Definition2.4} or Remark (7) of 
Section 2. If for some $\Gamma\subseteq\bigcap_{k\in {\Bbb N}} D(A_{n_k} 
)\cap {\cal C}\cap B_b(E)$ 
and for $g\in\Gamma$ the sequence $A_{n_k}g$ $s$-converges 
to zero as $k\to\infty$ then for $\t g\in {\cal C}$ we have    
\begin{eqnarray}\label{3.1}
\left\langle g,\t g\right\rangle_{n_k}&&\hspace{-.5cm}=\left\langle\beta 
G_{n_k,\beta}g,\t g\right\rangle_{n_k}-\left\langle G_{n_k,\beta}A_{n_k} 
g,\t g\right\rangle_{n_k} \nonumber \\ 
&&\hspace{-1cm}\stack{k\to\infty}{\lra}\int\beta G_\beta g\cdot\t g\, 
d\t\bnu\, , \quad\beta >0.  
\end{eqnarray}
%
This implies $g=\beta G_\beta g$, $\beta >0$, $\t\bnu$-a.e. and thus 
\begin{eqnarray}\label{3.2}
g=T_t g\, , \quad g\in\Gamma,\ t\ge 0,\ \t\bnu\mbox{\rm -a.e.}  
\end{eqnarray}
%
Let us recall that the trajectories $X_t$, $t\ge 0$, relative to $(T_t 
)_{t\ge 0}$ are supposed to be cadlag. Assuming even more that they are 
non-random given the initial value and that $\Gamma$ separates the points 
in $E$ (by containing an adequate subset of continuous functions) then 
(\ref{3.2}) implies $X_t=X_0$, $t\ge 0$, $\t\bnu$-a.e. 

Let the above hold for all accumulation points $\t\bnu$ of $\bnu_n$, 
$n\in {\Bbb N}$. 
\medskip 

If, as for example in the application of Subsection 4.1, the limiting 
process has just one stationary path taking the value $\mu_0$ then it 
follows that ${\bnu}_n\stack{n\to\infty}{\Ra}\delta_{\mu_0}=:\bnu$. 
Set $E_0:=E\setminus\bigcup_{n=1}^\infty E_n$. This introduction to 
the present section is enough motivation to introduce the following 
condition. 
\begin{itemize}
\item[(${\cal C}3$)] The measure ${\bnu}$ is concentrated on some $\mu_0 
\in E_0$. 
\end{itemize}

\subsection{Analysis on a Sequence of $L^2$-Spaces Continued}

In this subsection, we continue the work we have started in 
Proposition \ref{Proposition2.3}. We are interested in properties of 
$w$-convergent and $s$-convergent sequences which can be considered 
counterparts of properties of weak and strong convergent sequences 
in Hilbert spaces. Motivated by the above introduction to the present 
Section 3, throughout the whole Subsection 3.1, we will assume 
${\bnu}_n\stack{n\to\infty}{\Ra}{\bnu}$ for some probability measure 
$\bnu$ on $(E,{\cal B}(E))$. However, we do not necessarily assume 
$\bnu=\delta_{\mu_0}$. Let us specify the set ${\cal C}$ for this 
subsection. 
\begin{definition}\label{Definition3.1} 
A function $g\in\bigcap_{\, n\in {\Bbb N}}L^2(E,{\bnu}_n)\cap L^2 
(E,{\bnu})$ is said to {\it belong to the set} ${\cal C}$ if we have 
the following. 
\begin{itemize}
\item[(i)] There exsits a sequence $g_{0,r}\in C_b(E)$, $r\in {\Bbb N}$, 
with $g_{0,r}\stack{r\to\infty}{\lra}g$ in $L^2(E,\bnu)$ and $\limsup_{r 
\to\infty}\limsup_{n\to\infty}\langle g_{0,r},g_{0,r}\rangle_n$ $<\infty$. 
\item[(ii)] For every $n\in {\Bbb N}$, there exists a function $g_n \in 
C_b(E)$ such that $g=g_n$ ${\bnu}_n$-a.e., $\langle g_n,g_{0,r}\rangle_n 
\stack{n\to\infty}{\lra}\langle g,g_{0,r}\rangle$, $r\in {\Bbb N}$, and 
\item[(iii)] $\langle g_n,g_n\rangle_n\stack{n\to\infty}{\lra}\langle 
g,g\rangle$. 
\end{itemize} 
\end{definition} 
Obviously, ${\cal C}$ is a subset of ${\cal D}$ defined in Section 2. 
Linearity of ${\cal C}$ and thus the $s$-convergence of $g_n:=g\in 
{\cal C}$ to $g$ as $n\to\infty$ is left as an exercise. 
\begin{lemma}\label{Lemma3.2} 
Assume ${\bnu}_n\stack{n\to\infty}{\Ra}{\bnu}$ and ${\bnu}(\bigcup_{n 
\in {\Bbb N}}E_n)=0$. 
(a) We have $C_b(E)\subseteq {\cal C}$. In particular, condition 
(${\cal C}1$) is satisfied. \\ 
(b) Let 
\begin{eqnarray*}
f_k:=\left\{ 
\begin{array}{cl}
0 & \ \ \ \mbox{\rm on}\ \bigcup_{m=1}^{k-1}E_m\ \mbox{\rm for}\ 
k>1, \\ 
a_m^{(k)}\vp_k & \ \ \ \mbox{\rm on}\ E_m\, , \ m\ge k, 
\vphantom{\displaystyle \sum_1^1} \\ 
a^{(k)}\vp_k & \ \ \ \mbox{\rm on}\ E_0 
\end{array}
\right.\, , \quad k\in {\Bbb N}, 
\end{eqnarray*}
such that 
\begin{itemize}
\item[(i)] $\vp_k\in C_b(E)$, $k\in {\Bbb N}$, 
\item[(ii)] ${\Bbb R}\ni a_m^{(k)}\stack{m\to\infty}{\lra}a^{(k)}\in 
{\Bbb R}$, $k\in {\Bbb N}$, 
\item[(iii)] $C_1:=\sum_{k=1}^\infty\sup_{n\in {\Bbb N}}\langle f_k, 
f_k\rangle_n^{1/2}<\infty$. 
\end{itemize}
Then $g:=\sum_{k=1}^\infty f_k\in {\cal C}$. 
\end{lemma}
Proof. (a) This is an immediate consequence of ${\bnu}_n\stack{n\to 
\infty}{\Ra}{\bnu}$. 
\medskip

\nid 
(b) Let $g$ be defined as above and let $g_{0,r}:=\sum_{k=1}^ra^{(k)} 
\vp_k$. Taking into consideration ${\bnu}(\bigcup_{n\in {\Bbb N}}E_n) 
=0$, we verify condition (i) of Definition \ref{Definition3.1}. We shall 
demonstrate that $\langle g,g\rangle_n\stack{n\to\infty}{\lra}\langle g, 
g\rangle$. The proof of $\langle g,g_{0,r}\rangle_n\stack{n\to\infty} 
{\lra}\langle g,g_{0,r}\rangle$, $r\in {\Bbb N}$, is similar. 

\nid
{\it Step 1 } With $C_1$ defined in (iii), we have 
\begin{eqnarray*}
\left|\int\left(g^2-\left({\textstyle\sum_{k=1}^m} f_k\right)^2 
\right)\, d{\bnu}_n\right|&\le&2\sum_{i=1}^m\sum_{j=m+1}^\infty 
\int |f_if_j|\, d{\bnu}_n+\sum_{i,j=m+1}^\infty \int |f_if_j|\, 
d{\bnu}_n \nonumber \\ 
&&\hspace{-5.6cm}\le 2\sum_{i=1}^m\sum_{j=m+1}^\infty \left(\int 
f_i^2 \, d{\bnu}_n\right)^{1/2}\hspace{-2mm}\left(\int f_j^2\, 
d{\bnu}_n\right)^{1/2}\hspace{-1mm}+\hspace{-1mm}\sum_{i,j=m+1 
}^\infty \left(\int f_i^2\, d{\bnu}_n \right)^{1/2}\hspace{ 
-2mm}\left(\int f_j^2\, d{\bnu}_n\right)^{1/2}\hspace{-2mm} 
\nonumber \\ 
&&\hspace{-5.6cm}= 2\left(\sum_{i=1}^m \langle f_i\, , \, f_i 
\rangle_n^{1/2}\right) \left(\sum_{k=m+1}^\infty \langle f_{k} 
\, , \, f_{k}\rangle_n^{1/2}\right)+\left(\sum_{k=m+1}^\infty 
\langle f_{k}\, , \, f_{k}\rangle_n^{1/2}\right)^2 \nonumber \\ 
&&\hspace{-5.6cm}\le 2C_1\, \sum_{k=m+1}^\infty\langle f_{k}\, , \, 
f_{k}\rangle_n^{1/2}+\left(\sum_{k=m+1}^\infty\langle f_{k}\, , \, 
f_{k}\rangle_n^{1/2}\right)^2\, .  
\end{eqnarray*}
Now, this chain of equations as well as inequalities and (iii) show 
\begin{eqnarray*}
\left|\int\left(g^2-\left({\textstyle\sum_{k=1}^m} f_k\right)^2 
\right)\, d{\bnu}_n\right|&\stack{m\to\infty}{\lra}&0 \quad 
\mbox{uniformly in }\, n\in {\Bbb N}. 
\end{eqnarray*}
{\it Step 2 } On the other hand, with (ii), we observe that 
\begin{eqnarray*}
\int\left({\textstyle\sum_{k=1}^m} f_k\right)^2\, d{\bnu}_n&=& 
\sum_{i,j=1}^{m\wedge n}a_n^{(i)}a_n^{(j)}\int \vp_i\vp_j\, d{\bnu}_n\\ 
&\stack{n\to\infty}{\lra}& \sum_{i,j=1}^ma^{(i)}a^{(j)}\int 
\vp_i\vp_j\, d{\bnu}\\ 
&=&\int\left({\textstyle\sum_{k=1}^m}f_k\right)^2\, d{\bnu}\, . 
\end{eqnarray*}
{\it Step 3 } The results of Steps 1 and 2 finally yield 
\begin{eqnarray*}
\langle g\, , \, g\rangle_n=\int g^2\, d{\bnu}_n\stack{n\to\infty} 
{\lra}\int g^2\, d{\bnu}=\langle g\, , \, g\rangle\, . 
\end{eqnarray*}
\qed 
\begin{proposition}\label{Proposition3.3} 
Assume ${\bnu}_n\stack{n\to\infty}{\Ra}{\bnu}$. 
(a) Suppose ${\bnu}(\bigcup_{n\in {\Bbb N}}E_n)=0$. Let $\vp_n\in 
{\cal C}$, $n\in {\Bbb N}$, be a sequence $w$-convergent to $\vp\in 
L^2(E,{\bnu})$ as $n\to\infty$. Then $\langle \vp_n\, , \, \vp_n 
\rangle_n$, $n\in {\Bbb N}$, is bounded. \\ 
(b) Suppose again ${\bnu}(\bigcup_{n\in {\Bbb N}}E_n)=0$. Let $\vp_n 
\in {\cal C}$, $n\in {\Bbb N}$, be a sequence $w$-convergent to $\vp 
\in L^2(E,{\bnu})$ and let $\psi_n\in {\cal C}$, $n\in {\Bbb N}$, be a 
sequence that $s$-converges to $\psi\in L^2(E,\bnu)$ as $n\to\infty$. 
Then $\langle\vp_n\, , \, \psi_n\rangle_n\stack {n\to\infty}{\lra} 
\langle\vp\, , \, \psi\rangle$. \\ 
(c) Let $\t C_b(E)$ be a subset of $C_b(E)$ which is dense in $L^2(E, 
\bnu)$. For every $\vp\in L^2(E,\bnu)$, there is a subsequence $n_k$, 
$k\in {\Bbb N}$, of indices  and $\vp_{n_k}\in\t C_b(E)$, $k\in {\Bbb 
N}$, $s$-converging to $\vp$ and converging to $\vp$ in $L^2(E,\bnu)$ , 
both as $k\to\infty$. 
\end{proposition} 
Proof. (a) Without loss of generality, we may assume that $\langle\vp_n 
\, , \, \vp_n\rangle_m>0$ and $\langle\vp_n\, , \, \vp_n\rangle>0$ for 
all $m,n\in {\Bbb N}$. Otherwise, there is a suitable element of $C_b(E 
)$ which we add to all $\vp_n$. 

We shall show that $\langle \vp_n\, , \, \vp_n 
\rangle_n$, $n\in {\Bbb N}$, is bounded. However, let us assume that 
this was not the case. Then there is a subsequence $\vp_{n_k}$, $k\in 
{\Bbb N}$, of $\vp_n$, $n\in {\Bbb N}$, such that $\langle\vp_{n_k}\, , 
\, \vp_{n_k} \rangle_{n_k}^{1/2} \ge 4^k$, $k\in {\Bbb N}$. In the 
sense of the above definition of elements belonging to ${\cal C}$, let 
$\vp_{n_k}=g_{n_k}$ ${\bnu}_{n_k}$-a.e. where $g_{n_k}\in C_b(E)$. 
Since $\langle\vp_{n_k}\, , \, \vp_{n_k} \rangle_{n_k}^{1/2}=\langle 
g_{n_k}\, , \, g_{n_k} \rangle_{n_k}^{1/2}$, $k\in {\Bbb N}$, we may 
replace $\vp_{n_k}$ by $g_{n_k}$ in the sequence $\langle\vp_{n_k}\, , 
\, \vp_{n_k} \rangle_{n_k}^{1/2}$, $k\in {\Bbb N}$. Therefore, without 
loss of generality, we may assume that $\vp_{n_k}\in C_b(E)$, $k\in 
{\Bbb N}$. 
\medskip

Our goal is now to show that there is a $\psi\in {\cal C}$ such that 
$\langle\vp_{n_k}\, , \, \psi\rangle_{n_k}$ diverges as $k\to\infty$. For 
this introduce 
\begin{eqnarray*}
\psi_k:=\left\{ 
\begin{array}{cl}
0 & \ \ \ \mbox{\rm on}\ \bigcup_{m=1}^{k-1}E_m\ \mbox{for}\ k>1, 
\vphantom{{\displaystyle\frac{\vp_{n_k}}{\langle \vp_{n_k}\, , \, 
\vp_{n_k}\rangle^{1/2}_m}}}\\ 
\displaystyle\frac{\vp_{n_k}}{\langle \vp_{n_k}\, , \, \vp_{n_k} 
\rangle^{1/2}_{m}} & \ \ \ \mbox{\rm on}\ E_m\, , \ m\ge k, \\ 
\displaystyle\frac{\vp_{n_k}}{\langle \vp_{n_k}\, , \, \vp_{n_k} 
\rangle^{1/2}} & \ \ \ \mbox{\rm on}\ E_0 \vphantom{\displaystyle 
\sum^2}
\end{array}
\right.\, , \quad k\in {\Bbb N}, 
\end{eqnarray*}
and 
\begin{eqnarray*}
\psi=\sum_{k=1}^\infty \frac{b_k}{3^k}\psi_k
\end{eqnarray*}
where $b_k\in \{-1,1\}$ is chosen in such a way that 
\begin{eqnarray}\label{3.3}
b_k\sum_{r=1}^{k-1}\frac{b_r}{3^r} \langle \vp_{n_k}\, , \, \psi_r 
\rangle_{n_k} \ge 0\, , \quad k\in {\Bbb N}. 
\end{eqnarray}
Note that because of $\vp_{n_k}\in C_b(E)$, 
\begin{eqnarray*}
a_m^{(k)}:=\frac{b_k}{3^k\langle \vp_{n_k}\, , \, \vp_{n_k} 
\rangle^{1/2}_{m}}\stack{m\to\infty}{\lra}\frac{b_k}{3^k\langle 
\vp_{n_k}\, , \, \vp_{n_k}\rangle^{1/2}}=:a^{(k)}\, , \quad 
k\in {\Bbb N}, 
\end{eqnarray*}
and that with $f_k:=\frac{b_k}{3^k}\psi_k$, $k\in {\Bbb N}$, we have 
$\langle f_k\, , \, f_k\rangle_m^{1/2}\le\frac{1}{3^k}$ independent of 
$m\in {\Bbb N}$ and therefore 
\begin{eqnarray*}
\sum_{k=1}^\infty \sup_{m\in {\Bbb N}}\langle f_k\, , \, f_k\rangle_m^{ 
1/2}\le\sum_{k=1}^\infty\frac{1}{3^k}=\frac12.  
\end{eqnarray*}
With the above convention $\vp_{n_k}\in C_b(E)$, $k\in {\Bbb N}$, 
it follows now from Lemma \ref{Lemma3.2} (b) that $\psi\in {\cal C}$. 
Recalling the definition of $\psi$, relation (\ref{3.3}), and using 
$\langle \psi_r\, , \, \psi_r \rangle_{n_k}\in\{0,1\}$ as well as 
$\langle\vp_{n_k}\, , \, \psi_k\rangle_{n_k}=\langle\vp_{n_k}\, , \, 
\vp_{n_k} \rangle_{n_k}^{1/2}$, $r,k\in {\Bbb N}$, we obtain 
\begin{eqnarray*}
|\langle\vp_{n_k}\, , \, \psi\rangle_{n_k}|&=&\left|\sum_{r=1}^\infty
\frac{b_r}{3^r}\cdot \langle\vp_{n_k}\, , \, \psi_r\rangle_{n_k}\right| 
\\ 
&\ge&\frac{1}{3^k}\langle\vp_{n_k}\, , \, \vp_{n_k}\rangle_{n_k}^{1/2}- 
\sum_{r=k+1}^\infty\frac{1}{3^r}\left|\langle\vp_{n_k}\, , \, \psi_r 
\rangle_{n_k}\right|\\ 
&\ge&\frac{1}{3^k}\langle\vp_{n_k}\, , \, \vp_{n_k}\rangle_{n_k}^{1/2} 
-\sum_{r=k+1}^\infty\frac{1}{3^r}\langle\vp_{n_k}\, , \, \vp_{n_k} 
\rangle_{n_k}^{1/2} \\ 
&\ge&4^k\left(\frac{1}{3^k}-\sum_{r=k+1}^\infty\frac{1}{3^r}\right)\\ 
&=&\frac12 \left(\frac43\right)^k\, , \quad k\in {\Bbb N}.  
\end{eqnarray*}
Thus, $\langle\vp_{n_k}\, , \, \psi \rangle_{n_k}$ diverges as $k\to 
\infty$. Consequently, the assumption of Step 1 yields a contradiction, 
which in the end means that $\langle \vp_n\, , \, \vp_n\rangle_n$, 
$n\in {\Bbb N}$, is bounded. 
\medskip

\nid 
(b) Let $\ve >0$ and let $\t \psi\in {\cal C}$ with $\langle\psi -\t\psi 
\, , \, \psi -\t\psi\rangle^{1/2}<\ve$. Since $\psi_n\in {\cal C}$ is 
$s$-convergent to $\psi\in L^2(E,{\bnu})$ (which implies that $\psi_n$, 
$n\in{\Bbb N}$, is $w$-convergent to $\psi$), we obtain 
\begin{eqnarray*}
\langle\psi_n-\t \psi\, , \, \psi_n-\t \psi\rangle_n&=&\langle\psi_n\, , 
\, \psi_n\rangle_n-2\langle\psi_n\, , \, \t \psi\rangle_n+\langle\t \psi 
\, , \, \t \psi\rangle_n \\ 
&\stack{n\to\infty}{\lra}& \langle\psi -\t\psi\, , \, \psi -\t\psi\rangle 
\, .  
\end{eqnarray*}
Because of this relation and because of part (a) of this proposition, 
we get 
\begin{eqnarray*}
\limsup_{n\to\infty}|\langle\vp_n\, , \, \psi_n-\t \psi\rangle_n|&\le& 
\limsup_{n\to\infty}\left(\langle\vp_n\, , \, \vp_n\rangle^{1/2}_n\cdot 
\langle\psi_n-\t \psi\, , \, \psi_n-\t \psi\rangle^{1/2}_n\right) \\ 
&\le &\ve\cdot\sup_n\langle\vp_n\, , \, \vp_n\rangle^{1/2}_n\, ,  
\end{eqnarray*}
which implies 
\begin{eqnarray*}
&&\hspace{-1cm}\limsup_{n\to\infty}\left|\langle\vp_n\, , \, \psi_n\rangle_n 
-\langle\vp\, , \, \psi\rangle\right| \\ 
&&\le\lim_{n\to\infty}\left|\langle\vp_n\, , \, \t\psi\rangle_n-\langle 
\vp\, , \, \t\psi\rangle\right|+\left|\langle\vp\, , \, \psi-\t\psi\rangle 
\right|+\limsup_{n\to\infty}\left|\langle\vp_n\, , \, \psi_n-\t\psi\rangle_n 
\right| \\ 
&&\le\ve\cdot\left(\langle\vp\, , \, \vp\rangle^{1/2}+\sup_n\langle\vp_n\, 
, \, \vp_n\rangle^{1/2}_n\right)\, . 
\end{eqnarray*}
(c) Let $\t\vp_n\in \t C_b(E)$, $n\in {\Bbb N}$, such that $\langle\vp -\t 
\vp_n\, , \, \vp -\t\vp_n\rangle<\frac1n$, $n\in {\Bbb N}$. Let $\rho_i\in 
C_b(E)$, with $\|\rho_i\|=1$, $i\in {\Bbb N}$, be a sequence whose linear 
span is dense in $C_b(E)$. It holds that $\sum_{i=1}^\infty 2^{-i}\left| 
\langle\t\vp_n\, , \, \rho_i\rangle_m-\langle\t\vp_n\, , \, \rho_i\rangle 
\right|\stack{m\to\infty}{\lra}0$ and that $\langle\t\vp_n\, , \, \t\vp_n 
\rangle_m-\langle\t\vp_n\, , \, \t\vp_n\rangle\stack{m\to\infty}{\lra}0$, 
$n\in {\Bbb N}$. Thus, for $l,n\in {\Bbb N}$, there is an $n_l\ge n$ with 
$n_l>n_{l-1}$ if $l\ge 2$ such that we have the two inequalities $\sum_{i 
=1}^\infty 2^{-i}\left|\langle\t\vp_n\, , \, \rho_i\rangle_{n_l}-\langle 
\t\vp_n\, , \, \rho_i\rangle\right|<\frac1l$ and $\left|\langle\t\vp_n\, , 
\, \t\vp_n\rangle_{n_l}-\langle\t\vp_n\, , \, \t\vp_n\rangle\right|<\frac 
1l$. Setting $\vp_{n_l}:=\t\vp_n$, $l\in {\Bbb N}$, and keeping $\vp_{n_l} 
\stack{l\to\infty}{\lra}\vp$ in $L^2(E,\bnu)$ in mind, the second one of 
these two relations and Proposition \ref{Proposition2.3} (a) imply the 
existence of a subsequence $\vp_{n_k}$, $k\in {\Bbb N}$, of $\vp_{n_l}$, 
$l\in {\Bbb N}$, $w$-converging to some $\t\vp\in L^2(E,\bnu)$ but the 
first one says $\t\vp=\vp$. The latter yields the desired $s$-convergence. 
\qed
\medskip 

\nid
{\bf Remark} (Continuation of Remarks (7) and (8) of Section 2) We 
shall demonstrate how Proposition \ref{Proposition3.3} (b) may contribute 
to the verification of condition (i) of Definition \ref{Definition2.4}. 
For this assume ${\bnu}(\bigcup_{n\in {\Bbb N}}E_n)=0$ and that, for 
some set $\t C_b(E)\subseteq C_b(E)\bigcap_{n\in {\Bbb N}}D(A_n)$ dense 
in $L^2(E,\bnu)$, $(A,D(A))$ is the closure of $(A,\t C_b(E))$ on 
$L^2(E,\bnu)$. We will also suppose that $(A,D(A))$ is self-adjoint in 
$L^2(E,\bnu)$. Assume that in some application where we wish to verify 
condition (i) of Definition \ref{Definition2.4}, we have 
\begin{itemize}
\item[(ii'')] for $\psi\in\t C_b(E)$, $A_n\psi$ $s$-converges to $A\psi$ 
and $A'_n\psi$ $s$-converges to $A\psi$ as $n\to\infty$. 
\end{itemize}

Let $\vp$ and $\vp_n$, $n\in {\Bbb N}$, be as in condition (i) of 
Definition \ref{Definition2.4}. As already pointed out $(T_t)_{t\ge 0}$ 
is strongly continuous in $L^2(E,\bnu)$. Thus for $r,t>0$, 
\begin{eqnarray*} 
-\frac1t\int_{u=0}^tAT_u(rG_r\vp)\, du=\int_{s=0}^\infty re^{-rs}T_s 
\left(\frac{\vp-T_t\vp}{t}\right)\, ds\, . 
\end{eqnarray*} 
Using spectral representation, it follows that $\frac1t\int_{u=0}^tT_u 
(rG_r\vp)\, du\in D(A)$ and for $r,t>0$, 
\begin{eqnarray*} 
&&\hspace{-.5cm}-\int\left(2\vp-\frac1t\int_{u=0}^tT_u(rG_r\vp)\, du 
\right)\cdot A\left(\frac1t\int_{u=0}^tT_u(rG_r\vp)\, du\right)\, d 
\bnu\nonumber \\ 
&&\hspace{.5cm}=\int_{s=0}^\infty re^{-rs}\int\left(2\vp-\frac1t\int_{u 
=0}^tT_u(rG_r\vp)\, du\vphantom{\dot{f}}\right)\cdot T_s\left(\frac{\vp 
-T_t\vp}{t}\right)\, d\bnu\, ds \nonumber \\ 
&&\hspace{.5cm}=\int_{s=0}^\infty re^{-rs}\int_{\lambda=0}^\infty\left( 
2e^{-s\lambda}-\frac{e^{-s\lambda}}{t\lambda}\left(1-e^{-t\lambda}\right 
)\frac{r}{r+\lambda}\right)\cdot\frac{1-e^{-t\lambda}}{t}\, d\|E_\lambda 
\vp\|^2_{L^2(E,\sbnu)}\, ds\nonumber \\ 
&&\hspace{.5cm}=\int_{\lambda=0}^\infty\left(2\frac{r}{r+\lambda}-\frac{ 
1-e^{-t\lambda}}{t\lambda}\left(\frac{r}{r+\lambda}\right)^2\right)\cdot 
\frac{1-e^{-t\lambda}}{t}\, d\|E_\lambda\vp\|^2_{L^2(E,\sbnu)}
\end{eqnarray*} 
where $E_\lambda$, $\lambda\ge 0$, are the projection operators relative 
to the spectral resolution of $-A$. Now note that $\frac{1-e^{-t\lambda}} 
{t}\uparrow\lambda$ as $t\downarrow 0$. For $\vp\in D(S^c)$ the measure 
$\lambda\, d\|E_\lambda\vp\|^2_{L^2(E,\sbnu)}$ is finite and for $\vp\not 
\in D(S^c)$ the measure $\lambda\, d\|E_\lambda\vp\|^2_{L^2(E,\sbnu)}$ is 
infinite. Thus, 
\begin{eqnarray*} 
-\int\left(2\vp-\frac1t\int_{u=0}^tT_u(rG_r\vp)\, du\right)\cdot A\left( 
\frac1t\int_{u=0}^tT_u(rG_r\vp)\, du\right)\, d\bnu 
\end{eqnarray*} 
must tend to a finite value if $\vp\in D(S^c)$ and to $+\infty$ if $\vp 
\not\in D(S^c)$ as $t\to 0$ and then $r\to\infty$. In other words, there 
is a sequence $\t\vp_l\in D(A)$, $l\in {\Bbb N}$, with $\t\vp_l\stack{l 
\to\infty}{\lra}\vp$ in $L^2(E,\bnu)$ and 
\begin{eqnarray*} 
S\left(\t\vp_l,2\vp-\t\vp_l\right)\stack{l\to\infty}{\lra}S^c \left( 
\vp,\vp\right)\quad\mbox{\rm if}\quad\vp\in D(S^c) 
\end{eqnarray*} 
as well as 
\begin{eqnarray*} 
S\left(\t\vp_l,2\vp-\t\vp_l\right)\stack{l\to\infty}{\lra}\infty\quad 
\mbox{\rm if}\quad\vp\not\in D(S^c)\, . 
\end{eqnarray*} 
By the construction of $(A,D(A))$ as the closure of $(A,\t C_b(E))$ 
in $L^2(E,\bnu)$ we even may assume that $\t\vp_l\in\t C_b(E)$, $l\in 
{\Bbb N}$. Fix $l\in {\Bbb N}$ for a moment. By (ii''), the sequences 
$A_n\t\vp_l$ as well as $A'_n\t\vp_l$ $s$-converge to $A \t\vp_l$ and 
by hypothesis (i) $\vp_n-\t\vp_l$ $w$-converges to $\vp-\t\vp_l$, all 
as $n\to\infty$. Thus, Proposition \ref{Proposition3.3} (b) implies 
$\left\langle -A_n\t\vp_l,\vp_n-\t\vp_l\right\rangle_n\stack{n\to\infty} 
{\lra}S(\t\vp_l,\vp-\t\vp_l)$ as well as $\left\langle -A'_n\t\vp_l, 
\vp_n-\t\vp_l\right\rangle_n\stack{n\to\infty}{\lra}S(\t\vp_l,\vp-\t 
\vp_l)$. We get therefore from (ii'') 
\begin{eqnarray*} 
\liminf_{n\to\infty}S_n(\vp_n,\vp_n)&&\hspace{-.5cm}=\liminf_{n\to 
\infty}\left(\left\langle -A_n\t\vp_l,\t\vp_l\right\rangle_n+\left 
\langle -A_n\left(\vp_n-\t\vp_l\right),\vp_n-\t\vp_l\right\rangle_n 
\right) \\ 
&&+\liminf_{n\to\infty}\left\langle -A_n\t\vp_l,\vp_n-\t\vp_l\right 
\rangle_n+\liminf_{n\to\infty}\left\langle -A'_n\t\vp_l,\vp_n-\t\vp_l 
\right\rangle_n \\ 
&&\hspace{-.5cm}\ge S (\t\vp_l,\t\vp_l)\vphantom{\dot{f}}+2S(\t\vp_l, 
\vp-\t\vp_l)=S (\t\vp_l,2\vp-\t\vp_l)\, , \quad l\in {\Bbb N}. 
\end{eqnarray*} 
The desired relation, $\liminf_{n\to\infty}S_n(\vp_n,\vp_n)\ge S^c( 
\vp,\vp)$ if $\vp\in D(S^c)$ and $S_n(\vp_n,\vp_n)\stack{n\to\infty} 
{\lra}$ $\infty$ if $\vp\not\in D(S^c)$ follows now. Subsequences 
$n_k$, $k\in {\Bbb N}$, of indices can be handled similarly. 

\subsection{Relations to the Theory of K. Kuwae and T. Shioya} 

Targeting on applications to spectral geometry, in the recent paper 
\cite{KS03}, especially Subsection 2.2 therein, terms describing weak 
and strong convergence in a certain collection of Hilbert spaces have 
been introduced, as well. As in our setting, 
properties of weak and strong convergence in Hilbert spaces have 
been adapted for the development of their framework. However, 
there is a basic difference between their and our approach to weak 
convergence in such a collection of Hilbert spaces. 
\medskip 

Their definition (Definition 2.5 together with Definition 2.4 and Lemma 
2.1 (1), all in \cite{KS03}) would read in our context as follows: 
\medskip

\nid 
{\bf Definition} (a) A sequence $\psi_n\in L^2(E,{\bnu}_n)$, $n\in 
{\Bbb N}$, is said to be {\it s-convergent} to $\psi\in L^2(E,{\bnu})$ 
as $n\to \infty$ if 
\begin{itemize} 
\item[(i)] $\langle\psi_n-\psi \, , \, \psi_n-\psi\rangle_n\stack{n\to 
\infty}{\lra}0$. 
\end{itemize}
(b) A sequence $\vp_n\in L^2(E,{\bnu}_n)$, $n\in {\Bbb N}$, is said to 
be {\it w-convergent} to $\vp\in L^2(E,{\bnu})$ as $n\to\infty$ if 
\begin{itemize} 
\item[(i)] $\langle\vp_n\, , \, \psi_n\rangle_n\stack{n\to\infty}{\lra} 
\langle\vp\, , \, \psi\rangle$ for all $\psi\in L^2(E,{\bnu})$ and all 
sequences $\psi_n\in L^2(E,{\bnu}_n)$, $n\in {\Bbb N}$, $s$-converging to 
$\psi$. 
\end{itemize} 
This difference in the definition of weak convergence in a collection 
of Hilbert spaces has the following consequences. 
\begin{itemize} 
\item[(1)] Fixing a set ${\cal C}$ as in our setting as the set of test 
functions for $w$-convergence reflects the definition of weak convergence 
in Hilbert spaces in a natural way. Strictly speaking, we define 
{\it $w$-convergence relative to} ${\cal C}$ and consequently also 
{\it $s$-convergence relative to} ${\cal C}$, cf. Definition 
\ref{Definition2.2}. 
\item[(2)] Allowing the possibility to specify the set ${\cal C}$ 
according to the applications considered, see Definition 
\ref{Definition3.1} above, results in a specific framework of properties 
related to $w$- and $s$-convergence and in specific versions of the 
conditions of Mosco type convergence.   
\item[(3)] The price we have to pay for this is the fact that proving an 
analogy to their defining property (i) in the above Definition (b) is 
sophisticated work when starting with our Definition \ref{Definition2.2} 
(a). For this, recall our Lemma \ref{Lemma3.2} and our Proposition 
\ref{Proposition3.3}. For example, we cannot follow their idea to 
prove our Proposition \ref{Proposition3.3} (a). 
Attempting to adopt the proof of Lemma 2.3 in \cite{KS03} for this, we 
would necessarily need our Proposition \ref{Proposition3.3} (b) which is 
in our setting a consequence of Proposition \ref{Proposition3.3} (a) but 
in theirs the definition of $w$-convergence. 
\end{itemize} 

\subsection{Weak Convergence of Invariant Measures}

We are going to formulate consequences of Theorem \ref{Theorem2.7} under 
the condition (${\cal C}3$). Let the situation of Subsection 2.3 be in 
force. 
\begin{proposition}\label{Proposition3.4} 
(a) Let $\t {\bnu}$ be an accumulation point of $({\bnu}_n)_{n\in 
{\Bbb N}}$ and let $n_k$ be a subsequence of indices such that ${\bnu 
}_{n_k}\stack{k\to\infty}{\Ra}\t {\bnu}$. Let ${\cal C}$ be as introduced 
in Definition \ref{Definition3.1} with $\bnu$ replaced by $\t\bnu$. Let 
$\hat{C}(E)$ be a linear subset of ${\cal C}$ and $\hat{C}_b(E):=\hat{C} 
(E)\cap B_b(E)$ such that $\hat{C}_b(E)$ is dense in $L^2(E,\t\bnu)$. 
Assume the following. 
\begin{itemize}
\item[(i)] For every $g\in\hat{C}_b(E)$, $\beta>0$, and $n\in {\Bbb N}$, 
there exists a $v_n\in \hat{C}_b(E)$ such that $G_{n, \beta}g=v_n$ ${\bnu 
}_n$-a.e.  
\item[(ii)] For $g\in \hat{C}_b(E)$ and any $\beta >0$ there is a $v 
\in\hat{C}_b(E)$ such that $G_\beta g=v$ $\t{\bnu}$-a.e. 
\item[(iii)] Assume that $S_{n_k}$ converges to $S_{\t\sbnu}$ or $S^c_{ 
\t\sbnu}$ in the sense of Definition \ref{Definition2.4} or Remark (7) 
of Section 2 with ${\cal C}$ replaced by $\hat{C}(E)$. 
\item[(iv)] There is a set $\Gamma\subseteq\bigcap_{k\in {\Bbb N}} 
D(A_{n_k})\cap\hat{C}_b(E)$ 
that separates the points in $E$ such that $A_{n_k}g$ 
$s$-converges to zero as $k\to\infty$ for $g\in\Gamma$. 
\item[(v)] The trajectories $X_t$, $t\ge 0$, relative to $(T_t)_{t\ge 
0}$ regarded as a semigroup in $B_b(E)$ are non-random given the initial 
value and there is exactly one $\mu_0\in E$ such that $X_0=\mu_0$ implies 
$X_t=\mu_0$, $t>0$. 
\end{itemize}
Let (${\cal C}3$) and (v) be verified. For all accumulation points $\t 
\bnu$ of $({\bnu}_n)_{n\in {\Bbb N}}$ assume conditions (i)-(iv). Then 
${\bnu}_n\stack{n\to\infty}{\Ra}{\bnu}=\delta_{\mu_0}$. 
\medskip 

\nid 
(b) Assume ${\bnu}_n\stack{n\to\infty}{\Ra}{\bnu}=\delta_{\mu_0}$. Then 
we have 
\begin{eqnarray}\label{3.4}
\langle\beta G_{n,\beta}g-g\, , \, \beta G_{n,\beta}g-g\rangle_n 
\stack{n\to\infty}{\lra} 0\, , \quad g\in C_b(E). 
\end{eqnarray} 
\end{proposition}
Proof. (a) We have to adjust the ideas of the introduction to the present 
Section 3 to the situation in this part of the lemma. The crucial step is 
(\ref{3.1}). In particular, we have to take into consideration the choice 
of $\hat{C}_b(E)$. 
\medskip 

\nid
(b) It holds that 
\begin{eqnarray*}
&&\hspace{-.5cm}\langle \beta G_{n,\beta} g-g\, , \, \beta G_{n,\beta} 
g-g\rangle_n \vphantom{\sum} \\ 
&&\hspace{.5cm}\le 2\langle \beta G_{n,\beta} g-g(\mu_0)\1\, , \, \beta 
G_{n,\beta} g-g(\mu_0)\1\rangle_n+2\langle g(\mu_0)\1-g\, , \, g(\mu_0)\1 
-g\rangle_n \vphantom{\sum} \\ 
&&\hspace{.5cm}\le 4\langle g(\mu_0)\1-g\, , \, g(\mu_0)\1-g\rangle_n 
\vphantom{\sum}\stack{n\to\infty}{\lra} 0\, , \quad g\in C_b(E). 
\end{eqnarray*} 
\qed
\medskip 

\nid 
{\bf Remarks} (1) Assume that the semigroups $(T_{n,t})_{t\ge 0}$, 
$n\in {\Bbb N}$, and $(T_t)_{t\ge 0}$ are {\it Feller}, i. e., 
they leave the space $C_b$ invariant. If we choose $\hat{C}(E):= 
C_b(E)$ then conditions (i) and (ii) of Proposition \ref{Proposition3.4} 
(a) become trivial by Lemma \ref{Lemma3.2} (a). 
\medskip 

\nid 
(2) Part (b) of Proposition \ref{Proposition3.4} describes the situation 
when we already know that $\bnu=\delta_{\mu_0}$ is the limiting measure. 
However to prove that the limiting measure is $\bnu=\delta_{\mu_0}$, 
Proposition \ref{Proposition3.4} (a) may be useful as we will demonstrate 
in the application of Subsection 4.1. In particular, it has to be proved 
that the limiting measure is concentrated on one single point. 

\subsection{Weak Convergence of Particle Processes} 

Let ${\cal M}_1(\overline{D})$ denote the space of all probability 
measures on $(\overline{D},{\cal B}(\overline{D}))$ where $D$ is a 
bounded $d$-dimensional domain or, more general, a bounded 
$d$-dimensional Riemannian manifold for some $d\in {\Bbb N}$. In 
addition, let ${\cal M}_\partial(\overline{D})$ be the quotient set of 
${\cal M}_1(\overline{D})$ with respect to measures on the boundary 
$\partial D$, that is the set of all equivalence classes $\mu$ such 
that $\nu_1,\nu_2\in\mu$ implies $\nu_1|_D=\nu_2|_D$. Let both spaces 
${\cal M}_1(\overline{D})$ and ${\cal M}_\partial(\overline{D})$ be 
endowed with the Prohorov metric. 

In particular, we will assume that $E$ is one of the compact spaces 
${\cal M}_1(\overline{D})$ or ${\cal M}_\partial(\overline{D})$. 
Furtheremore, for $E={\cal M}_1(\overline{D})$ and $n\in {\Bbb N}$, 
let ${\cal M}^n_1 (\overline{D})$ be the set of all measures $\mu$ 
in $E$ of the form $\mu=\frac{1}{n}\sum_{i=1}^n\delta_{z_i}$ where 
$z_1,\ldots ,z_n\in\overline{D}$ and $\delta_z$ denotes the Dirac 
measure concentrated on $z$. We set $E_n:={\cal M}^n_1(\overline{D}) 
\setminus\bigcup_{k=1}^{n-1}{\cal M}^k_1(\overline{D})$, $n\in 
{\Bbb N}$. 

If $E={\cal M}_\partial (\overline{D})$ and $n\in {\Bbb N}$, let 
${\cal M}^n_\partial (\overline{D})$ be the set of all measures 
$\mu$ in $E$ of the form $\mu=\frac{1}{n}\sum_{i=1}^n\delta_{z_i}$ 
where $z_1,\ldots ,z_n\in\overline{D}$. We will write $E_n:={\cal 
M}^n_\partial(\overline{D})\setminus\bigcup_{k=1}^{n-1}{\cal 
M}^k_\partial (\overline{D})$. Also note that in the case of 
$E={\cal M}_\partial (\overline{D})$, we identify all points belonging 
to $\partial D$ with each other. 

From now on, we will assume that there are Markov processes associated 
with the semigroups $(T_{n,t})_{t\ge 0}$, $n\in {\Bbb N}$, and $(T_t 
)_{t\ge 0}$: For $n\in {\Bbb N}$, let $X^n=((X^n_t)_{t\ge 0},(P^n_\mu 
)_{\mu\in {\cal M}^n_1 (\overline{D})})$ or $X^n=((X^n_t)_{t\ge 0}, 
(P^n_\mu )_{\mu\in {\cal M}^n_\partial (\overline{D})})$ be a process 
corresponding to the semigroup $(T_{n,t})_{t\ge 0}$ which takes values 
in $E_n$. To ensure well-definiteness of the latter we always assume 
$P^n_\mu (X^n_t\in {\cal M}^n_1 (\overline{D})\setminus E_n)=0$, 
$\mu\in {\cal M}^n_1 (\overline{D})$, or $P^n_\mu (X^n_t\in {\cal 
M}^n_\partial (\overline{D})\setminus E_n)=0$, $\mu\in {\cal M}^n_\partial 
(\overline{D})$, respectively, for all $t>0$. 

Let $X=((X_t)_{t\ge 0},(P_\mu )_{\mu\in E})$ be a process associated 
with the semigroup $(T_t)_{t\ge 0}$ which takes values in some subset 
of $E$. Suppose that the paths of the processes $X^n$, $n\in {\Bbb N}$, 
are cadlag. Let $\|\cdot\|$ denote the supremum norm in $B_b(E)$. 
\medskip

For a given sequence $\ve_n>0$, $n\in {\Bbb N}$, introduce 
\begin{eqnarray*}
B:=\bigcup_{n=1}^\infty\left\{\mu \in E_{n}:\left|g(\mu )-\beta G_{ 
n,\beta}g(\mu )\right|\ge\ve_n\|g\|\right\}\, . 
\end{eqnarray*}
Let $\psi_n\in L^\infty (E,{\bnu}_n)$, $n\in {\Bbb N}$, be a sequence 
of nonnegative uniformly bounded functions, i.e., there exists $C>0$ 
such that $\psi_n\le C$ ${\bnu}_n$-a.e. on $E$, $n\in {\Bbb N}$. In 
addition, suppose $\int\psi_n\, d{\bnu}_n=1$, $n\in {\Bbb N}$. Define 
the measures $P_{\psi_n\sbnu_n}:=\int_E P^n_\mu\, \psi_n 
(\mu ){\bnu}_n(d \mu)$, $n\in{\Bbb N}$, and $P_{\sbnu}:=\int_E P_\mu 
\, {\bnu}(d\mu)$, and introduce the processes ${\bf X}^n=((X^n_t)_{ 
t\ge 0},P_{\psi_n\sbnu_n})$ and ${\bf X}=((X_t)_{t\ge 0},P_{\sbnu})$. 
Moreover, let $\, {\Bbb E}_{\psi_n}$ be the expectation corresponding 
to $P_{\psi_n\sbnu_n}$, $n\in {\Bbb N}$. Finally, let $\tau_{B^c} 
\equiv\tau_{B^c}^n$ denote the first exit time of ${\bf X}^{n}$ from 
the set $B^c\cap E_{n}$, $n\in {\Bbb N}$. Let $T>0$, set 
\begin{eqnarray*}
\gamma_n\equiv\gamma_n(g):=\sup_{s\in [0,T+1]}\left|g(X^n_s)-\beta 
G_{n,\beta}g(X^n_s)\right|\, , \quad n\in {\Bbb N}.  
\end{eqnarray*}
In order to prove relative compactness of the family of processes 
$g({\bf X}^n)=((g(X^n_t))_{t\ge 0})$, 
$n\in{\Bbb N}$, we need one more technical condition. In particular, we 
specify the sequence $\ve_n>0$, $n\in {\Bbb N}$. 
\begin{itemize}
\item[(${\cal C}4$)] There exists an algebra $\t C_b(E)\subseteq C_b 
(E)$ containing the constant functions and separating points in $E$ 
and there is a function $\xi : (0,\infty)\times\t C_b(E)\times {\Bbb 
N}\to (0,\infty)$ with $\lim_{n\to\infty}\xi (a_n,g,n)=0$ if $\lim_{ 
n\to \infty}a_n=0$ such that, for $g\in\t C_b(E)$ and $\ve_n:=\xi 
(\langle g-\beta G_{n,\beta}g\, ,\, g-\beta G_{n,\beta}g\rangle_{n 
}^{1/2}\, ,g\, ,n)$, we have 
\begin{eqnarray*}
\, {\Bbb E}_{\psi_n}\left(e^{-\beta\tau_{B^c}}\right)\le\ve_n\, ,\quad 
n\in {\Bbb N}. 
\end{eqnarray*}
\end{itemize}
\begin{theorem}\label{Theorem3.5} 
Let $\psi_n\in L^\infty(E,{\bnu}_n)$, $n\in {\Bbb N}$, be a 
sequence of nonnegative functions which are uniformly bounded by some 
constant $C>0$ and satisfy $\int\psi_n\, d{\bnu}_n=1$, $n\in {\Bbb N}$. 
In addition, suppose (\ref{3.4}) (cf. Proposition \ref{Proposition3.4}) 
and (${\cal C}4$) for some $\beta >0$. \\ 
(a) For $g\in \t C_b(E)$, the family of processes $g({\bf X}^n)= 
((g(X^n_t))_{t\ge 0})$, 
$n\in {\Bbb N}$, is relatively compact. \\ 
(b) The family of processes ${\bf X}^n=((X^n_t)_{t\ge 0})$, 
$n\in {\Bbb N}$, is relatively compact. 
\end{theorem}
Proof. (a) We will apply Theorem 3.8.6 of S. N. Ethier, T. Kurtz  
\cite{EK86}. For $n\in {\Bbb N}$ and $t\ge 0$, let ${\cal F}^n_t$ 
denote the $\sigma$-algebra generated by the family $(X^n_s)_{0\le 
s\le t}$. In Steps 1 and 2 below, we will keep $n\in {\Bbb N}$ fixed. 
In Step 3, we will then pass to the limit as $n\to\infty$. 
\medskip

\nid 
{\it Step 1 } Let $g\in\t C_b(E)$, $0<\delta <1$, $0\le t\le T$, 
$0\le u\le \delta$, and $\beta >0$. Since $\beta G_{n,\beta}g, 
\beta G_{n,\beta}g^2\in D(A_n)$ and $X^n$ is Markov, it follows 
from (1.5) of \cite{Dy65} that $\, {\Bbb E}_{\psi_n}[\vp(X^n_{t+u}) 
-\vp(X^n_t)|{\cal F}^n_t]=\, {\Bbb E}_{\psi_n}[\int_t^{t+u}A_n\vp 
(X^n_s)\, ds|{\cal F}^n_t]$ for $\vp=\beta G_{n,\beta}g$ and $\vp= 
\beta G_{n,\beta}g^2$ and therefore 
\begin{eqnarray}\label{3.5}
&&\hspace{-0.5cm} 
\, {\Bbb E}_{\psi_n}\left[\left(g(X^{n}_{t+u})-g(X^{n}_t)\right)^2 
|{\cal F}^{n}_t\right] \nonumber \vphantom{\sum_k}\\ 
&&=\, {\Bbb E}_{\psi_n}\left[\left.(g(X^{n}_{t+u}))^2-(g(X^{n}_t))^2 
|{\cal F}^{n}_t\right]-2g(X^{n}_t)\, {\Bbb E}_{\psi_n}\left[g(X^{
n}_{t+u})-g(X^{n}_t)\right|{\cal F}^{n}_t\right]\nonumber 
\vphantom{\sum_k} \\ 
&&\le 2\, {\Bbb E}_{\psi_n}\left[\left.\sup_{s\in [0,T+1]}\left| 
g^2(X^{n}_s)-\beta G_{{n},\beta}g^2(X^{n}_s)\right|\right| 
{\cal F}^{n}_t\right]\nonumber \\ 
&&\hspace{4.5cm}+\, {\Bbb E}_{\psi_n}\left[\sup_{r\in [0,T]}\left. 
\int_r^{r+\delta}\left|A_{n}(\beta G_{{n},\beta}g^2)(X^{n}_s)\right| 
\, ds\right|{\cal F}^{n}_t\right]\nonumber \\
&&\hphantom{\le}+ 4\|g\|\, {\Bbb E}_{\psi_n}\left[\left.\sup_{s\in 
[0,T+1]}\left|g(X^{n}_s)-\beta G_{{n},\beta}g(X^{n}_s) 
\right|\right|{\cal F}^{n}_t\right]\nonumber \\ 
&&\hspace{4.5cm}+2\|g\|\, {\Bbb E}_{\psi_n}\left[\sup_{r\in [0,T]} 
\left.\int_r^{r+\delta}\left|A_{n}(\beta G_{{n},\beta}g)(X^{n}_s) 
\right|\, ds\right|{\cal F}^{n}_t\right]\nonumber \\ 
&&= 2 \, {\Bbb E}_{\psi_n}\left[\left.\gamma_n(g^2)\right|{\cal F 
}^{n}_t\right]+\, {\Bbb E}_{\psi_n}\left[\sup_{r\in [0,T]}\left.\beta 
\int_r^{r+\delta}\left|g^2(X^{n}_s)-\beta G_{{n},\beta}g^2( X^{n}_s) 
\right|\, ds\right|{\cal F}^{n}_t\right] \nonumber \\ 
&&\hphantom{\le}+ 4\|g\|\, {\Bbb E}_{\psi_n}\left[\left.\gamma_n(g) 
\right|{\cal F}^{n}_t\right]+2\|g\| \, {\Bbb E}_{\psi_n}\left[\sup_{ 
r\in [0,T]}\left.\beta\int_r^{r+\delta}\left|g(X^{n}_s)-\beta G_{{n}, 
\beta}g(X^{n}_s)\right|\, ds\right|{\cal F}^{n}_t\right] 
\hspace{-6mm} \nonumber \\ 
&&\le 2 \, {\Bbb E}_{\psi_n}\left[\left.\gamma_n(g^2)\right|{\cal F 
}^{n}_t\right]+\beta\delta \, {\Bbb E}_{\psi_n}\left[\left.\sup_{s\in 
[0,T+1]}\left|g^2(X^{n}_s)-\beta G_{{n},\beta}g^2(X^{n}_s)\right| 
\right|{\cal F}^{n}_t\right] \nonumber \\ 
&&\hphantom{\le}+ 4\|g\|\, {\Bbb E}_{\psi_n}\left[\left. \gamma_n(g) 
\right|{\cal F}^{n}_t\right]+2\|g\|\beta \delta \, {\Bbb E}_{\psi_n} 
\left[\left.\sup_{s\in [0,T+1]}\left|g(X^{n}_s)-\beta G_{{n},\beta} 
g(X^{n}_s)\right|\right|{\cal F}^{n}_t\right] \nonumber \\ 
&&= (2+\beta\delta)\cdot \, {\Bbb E}_{\psi_n}\left[\left. 
\gamma_n(g^2)\right|{\cal F}^{n}_t\right]+2\|g\|(2+\beta\delta)
\cdot \, {\Bbb E}_{\psi_n}\left[\left.\gamma_n(g) \right|{\cal F 
}^{n}_t\right]\, . \hspace{1.5mm} 
\vphantom{\int^6}
\end{eqnarray}
{\it Step 2 } We have 
\begin{eqnarray}\label{3.6}
\, {\Bbb E}_{\psi_n}\gamma_n(g)&=&\, {\Bbb E}_{\psi_n}\left(\chi_{\{ 
\gamma_n\le\ve_n\|g\|\}}\gamma_n(g)\right)+\, {\Bbb E}_{\psi_n}\left( 
\chi_{\{\gamma_n>\ve_n\|g\|\}}\gamma_n (g)\right)\vphantom{\sum} 
\nonumber \\ 
&\le&\ve_n\|g\|+2\|g\|\, {\Bbb E}_{\psi_n}\left(\chi_{\{\gamma_n>\ve_n 
\|g\|\}}\right) \vphantom{\sum}\nonumber \\ 
&=&\ve_n\|g\|+2\|g\|e^{\beta(T+1)}\, {\Bbb E}_{\psi_n}\left(e^{-\beta( 
T+1)}\chi_{\{\gamma_n>\ve_n\|g\|\}}\right)\vphantom{\sum^1}\nonumber \\ 
&\le&\ve_n\|g\|+2\|g\|e^{\beta(T+1)}\, {\Bbb E}_{\psi_n}\left(e^{- 
\beta\tau_{B^c}}\chi_{\{\gamma_n>\ve_n\|g\|\}}\right)\vphantom{\sum^1_1} 
\nonumber \\ 
&\le&\ve_n\|g\|\cdot\left(1+2e^{\beta(T+1)}\right)\, , 
\end{eqnarray}
where the last line is justified by condition (${\cal C}4$). It 
follows now from relation (\ref{3.4}) and the definition of $\ve_n$ 
in (${\cal C}4$) that 
\begin{eqnarray}\label{3.7}
\, {\Bbb E}_{\psi_{n}}\gamma_{n}(g)\stack{n\to\infty}{\lra}0\, . 
\end{eqnarray}
Since $\t C_b(E)$ is assumed to be an algebra, from this it also 
follows that 
\begin{eqnarray}\label{3.8}
\, {\Bbb E}_{\psi_{n}}\gamma_{n}(g^2)\stack{n\to\infty}{\lra}0\, . 
\end{eqnarray}
{\it Step 3 } Setting 
\begin{eqnarray*}
\g_n(\delta):=(2+\beta\delta)\cdot\gamma_n(g^2)+2\|g\|(2+\beta \delta) 
\cdot\gamma_n(g)\, , \quad n\in {\Bbb N}, 
\end{eqnarray*}
and taking into consideration (\ref{3.5}), we observe that  
\begin{eqnarray*}
\, {\Bbb E}_{\psi_{n}}\left[\left(g(X^{n}_{t+u})-g(X^{n}_t)\right)^2 
|{\cal F}^{n}_t\right]\le \, {\Bbb E}_{\psi_{n}}\left[\left. \g_{n} 
(\delta)\right|{\cal F}^{n}_t\right]\, , \quad n\in {\Bbb N}, 
\end{eqnarray*}
and from (\ref{3.7}) and (\ref{3.8}), we obtain 
\begin{eqnarray*}
\, {\Bbb E}_{\psi_{n}}\g_{n}(\delta)\stack{n\to\infty}{\lra}0\, .
\end{eqnarray*}
Relative compactness of the family $g({\bf X}^n)$, $n\in {\Bbb N}$, 
follows now from  \cite{EK86}, Theorems 3.7.2 and 3.8.6, and Remark 
3.8.7. In particular, we note that $g({\bf X}^n)$, takes values in 
the compact interval $[\inf g,\sup g]$, $n\in {\Bbb N}$. \\ 
(b) By the Stone-Weierstrass Theorem, $\t C_b(E)$, is dense in 
$C_b(E)$. Furthermore, $E$ is compact. With these observations in 
mind, the claim follows from (a) and \cite{EK86}, Theorem 3.9.1. 
\qed
\begin{corollary}\label{Corollary3.6} 
Let ${\bnu}_n\stack{n\to\infty}{\Ra}{\bnu}=\delta_{\mu_0}$ (cf. 
Proposition \ref{Proposition3.4}), suppose condition (${\cal C}4$), 
and assume that the functions $\psi_n$, $n\in {\Bbb N}$, satisfy the 
conditions of Theorem \ref{Theorem3.5}. Then the processes ${\bf X 
}^n$ converge weakly to ${\bf X}$ as $n\to\infty$. 
\end{corollary}
The proof is an adaption to the proof of Theorem 6 in \cite{Lo05}. 
\qed

\section{Application}
\setcounter{equation}{0} 

In this section, we will apply Proposition \ref{Proposition3.4}, 
Theorem \ref{Theorem3.5}, and Corollary \ref{Corollary3.6} to a 
physically relevant situation. We will keep the notation of Section 
3 and \cite{Lo05} Section 2.  

Let $\Pi(n)$ denote the set of all permutations of the numbers $1, 
\ldots ,n$. For any permutation $\pi=(\pi(1),\ldots ,\pi(n))$, any 
$z_1,\ldots,z_n \in {\Bbb R}^d$, and $z=( z_1,\ldots ,z_n)$ introduce 
$z^\pi:=(z_{\pi(1)},\ldots ,z_{\pi(n)})$ and, for $A\in {\cal B}({\Bbb 
R}^{n\cdot d})$, set $A^\pi:=\{z^\pi:\,z\in A\}$. 
\medskip

Let the reader be reminded of the definitions of ${\cal M}_1(\overline{D 
})$ and ${\cal M}^n_1(\overline{D})$, $n\in {\Bbb N}$, in the beginning 
of Subsection 3.4. In addition, denote by ${\cal M}_1(D)$ the set of all 
probability measures on $(D,{\cal B}(D))$, let ${\cal M}^n_1(D)$ be the 
set of all measures $\mu=\frac{1}{n}\sum_{i=1}^n\delta_{z_i}$ where $z_1 
,\ldots ,z_n \in D$, $n\in {\Bbb N}$. 
\medskip

Let ${\cal B}^{\Pi(n)}(D^n):=\{A\in{\cal B}(D^n):\, A^\pi=A$ for all 
$\pi\in\Pi(n)\}$ and, for $A\in {\cal B}^{\Pi(n)}(D^n)$, set ${\cal A} 
\equiv {\cal A}(A):=\{\nu\in {\cal M}^n_1(D):\,\nu=\frac{1}{n}\sum_{i=1 
}^n\delta_{z_i},\ z=(z_1,\ldots ,z_n)\in A\}$. 

In the following we will use the notation $(h,\mu)=\int h\, d\mu$, 
$\mu\in {\cal M}_1(D)$, $h\in L^1(D,\mu)$. If $\mu$ has a density with 
respect to the Lebesgue measure then we will also write $(h,\rho)$ 
instead of $(h,\mu)$. 

\subsection{A Ginzburg-Landau Type Diffusion}

In this subsection, we will apply the results of Section 3 to a class 
of interacting diffusions on the circle as introduced in S. Lu \cite{Lu94}, 
\cite{Lu95}, S. Olla, S. R. S. Varadhan \cite{OV91}, K. Uchiyama \cite{Uc94}, 
and S. R. S. Varadhan \cite{Va91}. 
\medskip

\nid
Let $V\in C^1({\Bbb R})$ be an even function with $V\ge0$ and $\lim_{z\to 
\infty}V(z)=0$. Set $\Psi(z):=-zV'(z)$, $z\in {\Bbb R}$, and suppose $\Psi 
(z)\ge 0$, $z\in {\Bbb R}$. Assume that for some $\alpha>0$ and $\beta>0$
\begin{eqnarray}\label{4.1} 
\lim_{r\to\infty}r\cdot\Psi(r^{1+\alpha}z)=\beta|z|^{-\beta}\, , \quad z\in  
\left[-{\textstyle\frac12},{\textstyle\frac12}\right]\setminus\{0\}, 
\end{eqnarray}
uniformly on every compact subset of $\left[-{\textstyle\frac12},{\textstyle 
\frac12}\right]\setminus\{0\}$. It follows immediately that $\beta\in (0,1)$ 
and $\alpha=(1-\beta)/\beta$. We also suppose that the convergence 
(\ref{4.1}) holds in $L^1\left([t,\infty);z^{-1}\, dz\right)$ uniformly with 
respect to $t\in {\cal T}$ for every compact subset ${\cal T}\subset\left[-{ 
\textstyle\frac12},{\textstyle\frac12}\right]\setminus\{0\}$. This yields 
\begin{eqnarray}\label{4.2} 
\lim_{r\to\infty}r\cdot V(r^{1+\alpha}z)=|z|^{-\beta}\, , \quad z\in  
\left[-{\textstyle\frac12},{\textstyle\frac12}\right]\setminus\{0\}, 
\end{eqnarray}
where the convergence holds uniformly on every compact subset of $\left[ 
-{\textstyle\frac12},{\textstyle\frac12}\right]\setminus\{0\}$. In addition, 
we assume 
\begin{eqnarray}\label{4.3} 
V(z)\le z^{-\beta}\, ,\quad z\in (0,1)\quad\mbox{\rm and}\quad V(z)\ge z^{ 
-\beta}\, ,\quad z\in [1,\infty),  
\end{eqnarray}
which is equivalent to $\Psi\le\beta V$. 

Let $S$ be the circle of unit circumference, $(\Omega,{\cal F},P)$ be a 
probability space, and let $\beta_1,\beta_2,\ldots\, $ be a sequence of 
independent one-dimensional standard Brownian motions with state space $S$ 
on $(\Omega,{\cal F},P)$. Assume that for every $n\in {\Bbb N}$ we have an 
$S^n$-valued random element $x^n$ independent of $\beta_1,\beta_2,\ldots\, 
$ whose distribution under $P$ we denote by $\nu'_n$. 

In particular we will consider the measures $\nu_n$ on $(S^n,{\cal B} 
(S^n))$ which are defined by 
\begin{eqnarray}\label{4.4}
d\nu_n(x):=\frac{1}{Z_n}\exp\left[-\sum_{i,j=1}^nV(n^{1+\alpha}(x_i-x_j)) 
\right]\, dx_1\ldots dx_n
\end{eqnarray}
where $Z_n$ is a normalization constant. Let the $n$-particle process 
$(x_1(t),\ldots ,x_n(t))\equiv (x_1^n(t),$ $\ldots ,x_n^n(t))$ with state 
space $S^n$ that starts with $(x_1^n(0),\ldots ,x_n^n(0)):=(x_1^n,\ldots , 
x_n^n)$ follow the SDE 
\begin{eqnarray*}
dx_i(t)=-n^{1+\alpha}\sum_{j:j\neq i}V'(n^{1+\alpha}(x_i(t)-x_j(t)))\, 
dt + d\beta_i(t)\, ,\quad t\ge 0,\ i=1,\ldots ,n.  
\end{eqnarray*}
It is characterized by the closure $({\cal E}_n,D({\cal E}_n))$ on $L^2 
(S^n,\nu_n)$ of the positive symmetric bilinear form 
\begin{eqnarray*} 
{\cal E}_n(f,f):=\frac12\int_{S^n}\sum_{i=1}^n\left(\frac{\partial f} 
{\partial x_i}\right)^2\, d\nu_n\, , \quad f\in C^\infty (S^n). 
\end{eqnarray*}
The closure $({\cal E}_n,D({\cal E}_n))$ is a quasi-regular Dirichlet 
form which is associated with a strongly continuous contraction 
semigroup $(\hat{T}_{n,t})_{t\ge 0}$ on $L^2( S^n,\nu_n)$. The measure 
$\nu_n$ is invariant for the semigroup $(\hat{T}_{n,t})_{t\ge 0}$. 
The corresponding generator has the form 
\begin{eqnarray}\label{4.5} 
{\cal L}_n:=\frac12\sum_{i=1}^n\frac{\partial^2}{\partial x_i^2}-n^{1 
+\alpha}\sum_{j\neq i}^nV'(n^{1+\alpha}(x_i-x_j))\frac{\partial} 
{\partial x_i}\, .  
\end{eqnarray}
{\bf Remarks} (1) We note that the scaled process  
\begin{eqnarray*}
s:=n^{-2\alpha}t\, ,\quad y_i(s):=n^{\alpha}\cdot x_i\left(n^{-2 
\alpha}t\right)
\end{eqnarray*}
satisfies the SDE 
\begin{eqnarray*}
dy_i(s)=-n\sum_{j:j\neq i}V'(n(y_i(s)-y_j(s)))\, ds + dB_i(s)\, , 
\end{eqnarray*}
where $B_i(s):=n^\alpha\cdot\beta\left(n^{-2\alpha}t\right)$, $t\ge 
0$, $i=1,\ldots ,n$. This scaling establishes the formal connection 
to the papers \cite{Lu94}, \cite{Lu95}, \cite{OV91}, \cite{Uc94}, 
and \cite{Va91}. 
\medskip 

\nid 
(2) The major difference to the papers \cite{Lu94}, \cite{Lu95}, 
\cite{OV91}, \cite{Uc94}, and \cite{Va91} is that the function $V$ 
is no longer of compact support. Under the latter assumption, for 
large $n\in {\Bbb N}$, the analysis of $\left(V(n(x_i-x_j))\right 
)_{i,j=1,\ldots ,n}$ is carried out on some neighborhood of $\left 
\{(x_1,\ldots ,x_n)=(a,\ldots ,a):\right.$ $\left.a\in {\Bbb R} 
\right\}$. In contrast, (\ref{4.2}) assumes asymptotic behavior of 
$\left(V(n^{1+\alpha}(x_i-x_j))\right)_{i,j=1,\ldots ,n}$ for any 
argument $(x_1,\ldots ,x_n)\not\in\left\{(a,\ldots ,a):a\in {\Bbb 
R}\right\}$.  
\medskip  

Let $x^n=\left((x^n(t))_{t\ge 0},\hat{P}^n_\eta\right)$ where $x^n 
(t)=(x_1^n(t),\ldots ,x_n^n(t))$, $t\ge 0$, denotes the associated 
diffusion starting with an initial configuration $(x_1^n(0),\ldots 
,x_n^n(0)):=(\eta_1,\ldots ,\eta_n)\equiv\eta\in S^n$. Let 
$\hat{E}_\eta$ denote the expectation relative to $\hat{P}^n_\eta$. 
Set $\hat{P}_{f\nu_n}:=\int_{S^n}\hat{P}_\eta\cdot f(\eta)\, \nu_n 
({d\eta})$ and let $\hat{E}_{f}$ denote the expectation relative to 
$\hat{P}_{f\nu_n}$, $n\in {\Bbb N}$. 
\medskip 

For $m\in {\Bbb N}$, fix $b_m$ such that $\exp((\cdot)^{1/m})$ is 
convex on $[b_m,\infty)$. The first objective is to derive a PDE 
for the paths of the limiting process, cf. Proposition 
\ref{Proposition4.3} below. For this we assume that we start the 
processes $x^n$ with probability measures $d\nu'_n=f_n\, d\nu_n$, 
cf. (\ref{4.4}), such that all $f_n\in C^1(S^n)$ are symmetric in 
the $n$ entries and, for all odd natural numbers $m>m_0$ for some 
$m_0\in {\Bbb N}$ and all $n>n_0$ for some $n_0\in {\Bbb N}$ 
\begin{eqnarray}\label{4.6} 
e^{-(m_0-1)}\le f_n\le e^{nV(0)-b_{m_0}}\quad\mbox{\rm and}\quad\int 
f_n(\log f_n)^m\, d\nu_n\le (An)^m\, , 
\end{eqnarray}
for some $A>0$. We have the following property. For all odd natural 
numbers $m>m_0$ and $n>n_0$, let $f_n(t,\cdot)$, $t\ge 0$, be the 
solution to $\frac{d}{dt}\vp(t,x)={\cal L}_n\vp(t,x)$ with $\vp(0, 
\cdot)=f_n$. The function $H_m(t,n,f_n):=\int f_n(t,\cdot)(\log f_n 
(t,\cdot))^m\, d\nu_n$ is nonincreasing in $t\ge 0$. The proof of 
this is elementary. One takes the derivative of $H_m(t,n,f_n)$ with 
respect to $t$, uses the above PDE, the corresponding Dirichlet form 
representation, and the fact that the first part of (\ref{4.6}) 
implies $-(m_0-1)\le\log(f_n(t,\cdot))$ for all $t\ge 0$. 
\medskip

Let $D:=S$, note that $D=\overline{D}$. For $n\in {\Bbb N}$ introduce 
the measures ${\bnu}'_n$ on $(E_n,{\cal B}(E_n))$ by ${\bnu}'_n({\cal 
A}(A)):=\nu'_{n}(A)$, $A\in {\cal B}^{\Pi (n)}(D^n)$. Similarly, we 
introduce the measures $\bnu_n$, $n\in {\Bbb N}$. Let the process 
$X^n=\left((X^n_t)_{t\ge 0},(P^n_\mu)_{\mu\in {\cal M}_1^n(D)}\right)$ 
be defined by 
\begin{eqnarray*}
X^n_t:=\frac{1}{n}\sum_{i=1}^n\delta_{x^n_i(t)}\quad \mbox{\rm and} 
\quad P^n_\mu (X^n_{t}\in {\cal A}(A)):=\hat{P}^n_\eta(x^n(t)\in A)
\, , \quad t\ge 0, 
\end{eqnarray*}
$A\in {\cal B}^{\Pi (n)}(D^n)$, $\mu:=\frac1n\sum_{j=1}^n\delta_{ 
\eta_j}$, $\eta=(\eta_1,\ldots,\eta_n)\in D^n=S^n$. 
\medskip 

\nid 
{\bf Remark} (3) Since we are interested in the weak limits of the 
invariant measures $\bnu_n$ and the stationary versions of the 
processes $X^n$, $n\in {\Bbb N}$, hypothesis (\ref{4.6}) is not a 
restriction to us. In contrary, we even get asymptotic properties 
for a whole class of initial measures, namely those satisfying 
(\ref{4.6}). This condition implies also that our orientation should  
be the strategy of \cite{Va91} rather than the more general but also 
more sophisticated calculus of \cite{Uc94}. 
\begin{lemma}\label{Lemma4.1}  
Let $\bnu'_n$ be a sequence of probability measures on $E_n$, $n\in 
{\Bbb N}$, with (\ref{4.6}) and let $\t \bnu'$ be an arbitrary 
accumulation point of $\bnu'_n$, $n\in {\Bbb N}$. Assume ${\bnu}'_{ 
n_k}\stack{k\to\infty}{\Ra}\t\bnu'$. \\ 
(a) There exists $B>0$ such that 
\begin{eqnarray*}
\sup_{t\in [0,1]}\sup_{k\in {\Bbb N}}\frac{1}{n_k}\sum_{i,j=1}^{n_k}V 
\left(n_k^{1+\alpha}(x_i^{n_k}(t)-x_j^{n_k}(t))\right)<B\quad P\mbox{ 
\rm -a.s.} 
\end{eqnarray*}
(b) The measure $\t\bnu'$ is concentrated on the set of all probability 
measures $\mu (d\theta)=\rho(\theta)\, d\theta$ on $(S,{\cal B}(S))$ 
satisfying 
\begin{eqnarray}\label{4.7} 
\int_{\theta\in S}\int_{\tau\in S}|\theta-\tau|^{-\beta}\, \mu(d\tau) 
\, \mu(d\theta)\le B\, . 
\end{eqnarray}
\end{lemma} 
Proof. (a) Let $f_n(t,\cdot)$, $t\ge 0$, be the solution to $\frac{d} 
{dt}\vp(t,x)={\cal L}_n\vp(t,x)$ with $\vp(0,\cdot)=f_n$ and let $m> 
m_0$ be an odd natural number. For all $k\in {\Bbb N}$ such that $n_k 
>n_0$ we have as a consequence of (\ref{4.6}), $n_kV(0)-\sup_{x}\log 
f_{n_k}(t,x)\ge b_m$. We obtain 
\begin{eqnarray*}
&&\hspace{-.5cm}\int\left(\frac{1}{n_k}\sum_{i,j=1}^{n_k}V\left( 
n_k^{1+\alpha}(x_i-x_j)\right)\right)^mf_{n_k}(t,x)\, \nu_{n_k}(dx) 
 \\ 
&&\hspace{.5cm}\le\frac{2^{m-1}}{n_k^m}\int\left(-\log f_{n_k}(t,x)+ 
\sum_{i,j=1}^{n_k}V\left(n_k^{1+\alpha}(x_i-x_j)\right)\right)^mf_{ 
n_k}(t,x)\, \nu_{n_k}(dx) \\ 
&&\hspace{1.0cm}+2^{m-1}\, \frac{H_m(t,n_k,f_{n_k})}{n_k^m}\nonumber 
 \\ 
&&\hspace{.5cm}=\frac{2^{m-1}}{n_k^m}(\log(\cdot))^m\circ\exp((\cdot 
)^{1/m})\left(\int\left(-\log f_{n_k}(t,x)+\sum_{i,j=1}^{n_k}V(n^{1+ 
\alpha}(x_i-x_j))\right)^m\times\right.  \\ 
&&\hspace{1.0cm}\left.\times f_{n_k}(t,x)\, \nu_{n_k}(dx)\vphantom{ 
\sum_{i,j=1}^{n_k}}\right)+2^{m-1}\, \frac{H_m(t,n_k,f_{n_k})}{n_k^m 
}\, . 
\end{eqnarray*} 
We recall the below (\ref{4.6}) mentioned property and (\ref{4.2}) to 
verify 
\begin{eqnarray}\label{4.8}
&&\hspace{-.5cm}2\, \int\left(\frac{1}{n_k}\sum_{i,j=1}^{n_k}V\left( 
n_k^{1+\alpha}(x_i-x_j)\right)\right)^mf_{n_k}(t,x)\, \nu_{n_k}(dx) 
\nonumber \\ 
&&\hspace{.5cm}\le\left(\frac{2}{n_k}\log\int\frac{1}{f_{n_k}(t,x)} 
\exp\left\{\sum_{i,j=1}^{n_k}V\left(n_k^{1+\alpha}(x_i-x_j)\right) 
\right\}f_{n_k}(t,x)\, \nu_{n_k}(dx)\right)^m+(2A)^m\nonumber \\ 
&&\hspace{.5cm}=\left(-\frac{2}{n_k}\log\int_S\ldots\int_S\exp\left\{ 
-\sum_{i,j=1}^{n_k}V\left(n_k^{1+\alpha}(x_i-x_j)\right)\right\}\, d 
x_1\ldots\, dx_{n_k}\right)^m+(2A)^m\nonumber \\ 
&&\hspace{.5cm}\le\left(2\, \frac{n_k-1}{n_k}\int_S\ldots\int_Sn_kV 
\left(n_k^{1+\alpha}(x_1-x_2)\right)\, dx_1\ldots\, dx_{n_k}+2V(0) 
\right)^m+(2A)^m\nonumber \\ 
&&\hspace{.0cm}\stack{k\to\infty}{\lra}(2A_1)^m+(2A)^m\vphantom 
{\int}
\end{eqnarray} 
where we have applied Jensen's inequality two times. Recalling that 
$x^{n_k}(t)$ is an $S^{n_k}$-valued random element whose distribution 
under $P$ is $f_{n_k}(t,x)\, \nu_{n_k}(dx)$, relation (\ref{4.8}) can 
also be written as 
\begin{eqnarray*}
2\, \limsup_{k\to\infty}E\left(\frac{1}{n_k}\sum_{i,j=1}^{n_k}V\left( 
n_k^{1+\alpha}(x^{n_k}_i(t)-x^{n_k}_j(t))\right)\right)^m\le (2A_1)^m 
+(2A)^m\, ,\quad m>m_0, 
\end{eqnarray*} 
where $E$ is the expectation with respect to $P$. In other words, 
there exists $B>0$ such that for all $t\in [0,1]$
\begin{eqnarray}\label{4.9}
\limsup_{k\to\infty}\frac{1}{n_k}\sum_{i,j=1}^{n_k}V\left(n_k^{1+ 
\alpha}(x_i^{n_k}(t)-x_j^{n_k}(t))\right)<B\quad P\mbox{\rm -a.s.} 
\end{eqnarray} 
The claim is now a standard consequence using the fact that 
$\limsup_{k\to\infty}F_k$ of continuous functions $F_k$ is the 
decreasing limit of the lower semi-continuous functions $\sup_{k 
\ge N}F_k$ as $N\to\infty$. 
\medskip 

\nid
(b) Let $\mu_{n_k}:=\frac{1}{n_k}\sum_{i=1}^{n_k}\delta_{x^{n_k 
}_i}$, $k\in {\Bbb N}$. For $t=0$ we obtain from (\ref{4.9})
\begin{eqnarray*}
&&\hspace{-.5cm}\limsup_{k\to\infty}\bnu'_{n_k}\left(\int_{\theta 
\in S}\int_{\tau\in S}n_kV\left(n_k^{1+\alpha}(\theta-\tau)\right) 
\, \mu_{n_k}(d\tau)\, \mu_{n_k}(d\theta)>B\right) \\  
&&\hspace{.5cm}\le E\left(\limsup_{k\to\infty}\1_{\left\{\int_{ 
\theta\in S}\int_{\tau\in S}n_kV\left(n_k^{1+\alpha}(\theta-\tau) 
\right)\, \mu_{n_k}(d\tau)\, \mu_{n_k}(d\theta)>B\right\}}\right) 
 \\  
&&\hspace{.5cm}=E\left(\1_{\left\{\limsup_{k\to\infty}\int_{\theta 
\in S}\int_{\tau\in S}n_kV\left(n_k^{1+\alpha}(\theta-\tau)\right) 
\, \mu_{n_k}(d\tau)\, \mu_{n_k}(d\theta)>B\right\}}\right)=0\, . 
\end{eqnarray*} 
Let us define
\begin{eqnarray*} 
{\cal N}_B^k:=\left\{\mu_{n_k}\in E_{n_k}:\int_{\theta\in S}\int_{ 
\tau\in S}n_kV\left(n_k^{1+\alpha}(\theta-\tau)\right)\, \mu_{n_k} 
(d\tau)\, \mu_{n_k}(d\theta)\le B\right\}
\end{eqnarray*} 
and 
\begin{eqnarray*} 
{\cal N}_B:=\left\{\mu\in E:\int_{\theta\in S}\int_{\tau\in S}|\theta 
-\tau|^{-\beta}\, \mu(d\tau)\, \mu(d\theta)\le B\right\}\, . 
\end{eqnarray*} 
It is a straight consequence of the definitions of $\nu_{n_k}$, $\nu'_{ 
n_k}$, and $\bnu'_{n_k}$ that for the weak limit $\t \bnu'$, $\t \bnu'( 
\bigcup_{k\in {\Bbb N}}E_{n_k})=0$. By (\ref{4.2}), we have ${\cal N}_B 
\supseteq\overline{\bigcup_{k\in {\Bbb N}}{\cal N}_B^k}\setminus 
\bigcup_{k\in {\Bbb N}}E_{n_k}$, where $\overline{\bigcup_{k\in { 
\Bbb N}}{\cal N}_B^k}$ denotes the closure with respect to the weak 
topology in ${\cal M}_1(S)=E$. Here we have also taken advantage of 
the assumption that the convergence (\ref{4.2}) holds uniformly on 
every compact subset of $\left[-{\textstyle\frac12},{\textstyle 
\frac12}\right]\setminus\{0\}$. Therefore 
\begin{eqnarray*} 
\t \bnu'({\cal N}_B)\ge\t \bnu\left(\overline{\textstyle\bigcup_{k\in  
{\Bbb N}}{\cal N}_B^k}\right)\ge\liminf_{k\to\infty}\bnu'_{n_k}\left( 
\overline{\textstyle\bigcup_{k\in {\Bbb N}}{\cal N}_B^k}\right)\ge 
\liminf_{k\to\infty}\bnu'_{n_k}\left({\cal N}_B^k\right)=1
\end{eqnarray*} 
which completes the proof of part (b). 
\qed
\bigskip 

For the next lemma introduce the notation $E_{\rm abs}:=\{\mu\in E: 
\mu(d\theta)=\rho(\theta)\, d\theta, \rho\in L^1(S), \rho\ge 0\}$. 
\begin{lemma}\label{Lemma4.2} 
Let $\bnu'_n$ be a sequence of probability measures on $E_n$, $n\in 
{\Bbb N}$, with (\ref{4.6}) and let $\t \bnu'$ be an arbitrary 
accumulation point of $\bnu'_n$, $n\in {\Bbb N}$. Assume ${\bnu}'_{ 
n_k}\stack{k\to\infty}{\Ra}\t\bnu'$. Let $f\in C^1(S)$. \\ 
(a) The function $\Phi$ of type $\bigcup_{k\in {\Bbb N}}E_{n_k}\cup 
E_{\rm abs}\to [0,\infty]$ which is for $k\in {\Bbb N}$ and $\mu\in 
E_{n_k}$ defined by 
\begin{eqnarray*}
\Phi(\mu):=\int_{(\theta,\tau)\in S\times S\setminus D}\frac{f(\theta 
)-f(\tau)}{\theta -\tau}\cdot n_k\Psi\left(n_k^{1+\alpha}\left(\theta 
-\tau\right)\right)\, \mu(d\theta)\, \mu(d\tau) 
\end{eqnarray*}
and for $\mu(d\theta)=\rho(\theta)\, d\theta\in E_{\rm abs}$ by 
\begin{eqnarray*}
\Phi(\mu):=\beta\int_{\theta\in S}\int_{\tau\in S}\frac{f(\theta)-f 
(\tau)}{\theta -\tau}|\theta-\tau|^{-\beta}\rho (\theta)\rho(\tau)\, 
d\tau\, d\theta 
\end{eqnarray*}
is continuous with respect to the topology of weak convergence in ${ 
\cal M}_1(S)$ in every point $\mu(d\theta)=\rho(\theta)\, d\theta\in 
E_{\rm abs}$ for which $\Phi(\mu)<\infty$. \\ 
(b) We have 
\begin{eqnarray*}
\t \bnu'\left(\mu\in E:\Phi\ \mbox{\rm is discontinuous in }\mu\vphantom 
{l^1}\right)=0\, . 
\end{eqnarray*}
\end{lemma}
Proof. For part (a) recall that the convergence (\ref{4.1}) holds 
uniformly on every compact subset of $\left[-{\textstyle\frac12}, 
{\textstyle\frac12}\right]\setminus\{0\}$. Part (b) follows from part 
(a) of the present lemma and Lemma \ref{Lemma4.1} (b). 
\qed
\begin{proposition}\label{Proposition4.3}  
(a) The distribution $Q(\rho,\cdot)$ which is for all probability 
measures $\mu (d\theta)=\rho(\theta)\, d\theta$ on $(S,{\cal B}(S) 
)$ satisfying (\ref{4.7}) and all test functions $f\in C^1(S)$ given 
by 
\begin{eqnarray*}
Q(\rho,\cdot)(f)=\frac{\beta}{2}\int_{\theta\in S}\int_{\tau\in S} 
\frac{f(\theta)-f(\tau)}{\theta -\tau}|\theta-\tau|^{-\beta}\rho ( 
\theta)\rho (\tau)\, d\tau\, d\theta 
\end{eqnarray*}
has for all $\rho\in C^1(S)$ the representation 
\begin{eqnarray*}
Q(\rho,\cdot)(f)=\int_{\theta\in S}f(\theta)Q(\rho,\theta)\, d 
\theta
\end{eqnarray*}
where  
\begin{eqnarray*}
Q(\rho,\theta)=-\beta\int_{\tau\in (0,\frac12]}\frac{\rho(\theta+ 
\tau)-\rho(\theta-\tau)}{\tau}\, \tau^{-\beta}\, d\tau\cdot\rho( 
\theta)\, ,\quad\theta\in S. 
\end{eqnarray*}
(b) Let $\bnu'_n$ be a sequence of probability measures on $E_n$, $n 
\in {\Bbb N}$, with (\ref{4.6}) and let $\t \bnu'$ be an arbitrary 
accumulation point of $\bnu'_n$, $n\in {\Bbb N}$. For $\t \bnu'$-a.e. 
$\mu(d\theta)=\rho_0(\theta)\, d\theta$ the equation 
\begin{eqnarray}\label{4.10} 
\left\{
\begin{array}{rcl}
\D\int\frac{\partial}{\partial t}\rho(t,\theta)h(\theta)\, d\theta 
&=&\D\frac12\int h''(\theta)\rho(t,\theta)\, d\theta+Q(\rho(t,\cdot) 
,\cdot)(h') \\ 
\rho(t,\cdot)|_{t=0}&=&\rho_0\,  \vphantom{\D\int^1}
\end{array}\right.\, ,\quad h\in C^2(S), 
\end{eqnarray}
has a unique solution $\rho(t,\cdot)$, $t\ge 0$. For $\t \bnu'$-a.e. 
$\mu(d\theta)=\rho_0(\theta)\, d\theta$ and all $t\ge 0$, the measure 
$\mu(t,d\theta):=\rho(t,\theta)\, d\theta$ satisfies (\ref{4.7}). \\ 
(c) For $\t\bnu'$-a.e. $\mu(d\theta)=\rho(\theta)\, d\theta$ and fixed 
$t\ge 0$, the measure $\rho(t,\theta)\, d\theta$ depends continuously 
on the initial value $\rho(\theta)\, d\theta$ in the following sense. 

For each $\delta>0$ and $\rho(\theta)\, d\theta =\mu (d\theta)\in 
{\rm supp}\t\bnu'$ there is an $\ve >0$ such that if $\rho'(\theta) 
\, d\theta$ belongs to ${\rm supp}\t\bnu'$ and the $\ve$-neighborhood 
of $\rho(\theta)\, d\theta$ with respect to the Prohorov topology 
then $\rho'(t,\theta)\, d\theta$ belongs to the $\delta$-neighborhood 
of $\rho(t,\theta)\, d\theta$. \\ 
(d) Let $\mu_0$ be the measure on $(S,{\cal B}(S))$ which is the 
uniform distribution on $S$. Then the relation 
\begin{eqnarray}\label{4.11}  
0=\frac12\int h''(\theta)\rho(t,\theta)\, d\theta+Q(\rho(t,\cdot), 
\cdot)(h')\, ,\quad h\in C^2(S),\ t\ge 0,  
\end{eqnarray} 
implies $\rho(0,\theta)\, d\theta=\rho_0(\theta)\, d\theta=\mu_0 
(d\theta)$ with $\rho_0(\theta)=1$ for all $\theta\in\left[-\frac 
12,\frac12\right]$. 
\end{proposition} 
Proof. (a) For $f\in C^1(S)$ we have 
\begin{eqnarray*}
&&\hspace{-.5cm}\frac{\beta}{2}\int_{\theta\in S}\int_{\tau\in S} 
\frac{f(\theta)-f(\tau)}{\theta-\tau}|\theta-\tau|^{-\beta}\rho( 
\theta)\rho(\tau)\, d\tau\, d\theta \\ 
&&\hspace{.5cm}=\frac{\beta}{2}\lim_{\ve\to 0}\int_{\theta\in S} 
\int_{\tau:|\theta-\tau|\ge\ve}\frac{f(\theta)-f(\tau)}{\theta-\tau} 
|\theta-\tau|^{-\beta}\rho(\tau)\, d\tau\rho(\theta)\, d\theta \\ 
&&\hspace{.5cm}=\beta\lim_{\ve\to 0}\int_{\theta\in S}\int_{\tau:| 
\theta -\tau|\ge\ve}\frac{f(\theta)}{\theta-\tau}|\theta-\tau|^{- 
\beta}\rho(\tau)\, d\tau\rho(\theta)\, d\theta \\ 
&&\hspace{.5cm}=\beta\lim_{\ve\to 0}\int_{\theta\in S}f(\theta) 
\int_{\tau:|\theta -\tau|\ge\ve}\frac{|\theta-\tau|^{-\beta}}{\theta 
-\tau}\rho(\tau)\, d\tau\cdot\rho(\theta)\, d\theta \\ 
&&\hspace{.5cm}=-\beta\int_{\theta\in S}f(\theta)\left(\int_{\tau\in 
(0,\frac12]}\frac{\rho(\theta+\tau)-\rho(\theta-\tau)}{\tau}\, 
\tau^{-\beta}\, d\tau\cdot\rho(\theta)\right)\, d\theta \\ 
&&\hspace{.5cm}=\int_{\theta\in S}f(\theta)Q(\rho,\theta)\, d\theta 
\, . 
\end{eqnarray*}
(b) {\it Step 1 } Denote, more suggestively, 
\begin{eqnarray*} 
\mu_{n_k}(t,\cdot)\equiv X^{n_k}_t:=\frac{1}{n_k}\sum_{i=1}^{n_k} 
\delta_{x^{n_k}_i(t)}\, , \quad t\ge 0,\ k\in {\Bbb N}. 
\end{eqnarray*}

Assume ${\bnu}'_{n_k}\stack{k\to\infty}{\Ra}\t\bnu'$. Let $h\in 
C^2(S)$. For $k\in {\Bbb N}$ we obtain from It\^o's formula that  
$P$-a.s.
\begin{eqnarray}\label{4.12}
&&\hspace{-1.cm}\left(h,\mu_{n_k}(t,\cdot)\vphantom{l^1}\right)- 
\left(h,\mu_{n_k}(0,\cdot)\vphantom{l^1}\right)-\frac12\int_{s=0 
}^t\left(h'',\mu_{n_k}(s,\cdot)\vphantom{l^1}\right)\, ds-\int_{ 
s=0}^t\frac{1}{n_k}\sum_{i=1}^{n_k}h'\left(x^{n_k}_i(s)\right)\, 
d\beta_i(s)\nonumber \\ 
&&\hspace{.0cm}=-\int_{s=0}^tn_k^\alpha\sum_{i,j=1}^{n_k}h'\left 
(x^{n_k}_i(s)\right)\cdot V'\left(n_k^{1+\alpha}\left(x^{n_k}_i(s) 
-x^{n_k}_j(s)\right)\right)\, ds\nonumber \\ 
&&\hspace{.0cm}=-\int_{s=0}^t\frac{n_k^\alpha}{2}\sum_{i\neq j} 
\left(h'\left(x^{n_k}_i(s)\right)-h'\left(x^{n_k}_j(s)\right) 
\right)\cdot V'\left(n_k^{1+\alpha}\left(x^{n_k}_i(s)-x^{n_k}_j 
(s)\right)\right)\, ds\nonumber \\ 
&&\hspace{.0cm}=\int_{s=0}^t\frac{1}{2n_k}\sum_{i\neq j}\frac{h' 
\left(x^{n_k}_i(s)\right)-h'\left(x^{n_k}_j(s)\right)}{x^{n_k}_i 
(s)-x^{n_k}_j(s)}\cdot\Psi\left(n_k^{1+\alpha}\left(x^{n_k}_i(s) 
-x^{n_k}_j(s)\right)\right)\, ds\nonumber \\ 
&&\hspace{.0cm}=\int_{s=0}^t\int_{(\theta,\tau)\in S\times S 
\setminus D}\frac{h'(\theta)-h'(\tau)}{\theta-\tau}\cdot\frac{n_k 
}{2}\Psi\left(n_k^{1+\alpha}\left(\theta -\tau\right)\right)\, 
\mu_{n_k}(s,d\theta)\, \mu_{n_k}(s,d\tau)\, ds 
\end{eqnarray}
where, for the second equality sign, we have taken into 
consideration that $V'$ is skew symmetric with $V'(0)=0$. Let us 
take a closer look at the items of the first and the last line of 
(\ref{4.12}) which reads now as $\left(h,\mu_{n_k}(t,\cdot) 
\vphantom{l^1}\right)-\left(h,\mu_{n_k}(0,\cdot)\vphantom{l^1} 
\right)-I^{(1)}_{n_k}(t)-I^{(2)}_{n_k}(t)=I^{(3)}_{n_k}(t)$. 

Clearly, $\left(h,\mu_{n_k}(t,\cdot)\vphantom{l^1}\right)-\left( 
h,\mu_{n_k}(0,\cdot)\vphantom{l^1}\right)$ is $P$-a.e. uniformly 
bounded with respect to $t\in [0,1]$ and $k\in {\Bbb N}$. The 
expression $I^{(1)}_{n_k}(t)$ is $P$-a.e. equicontinuous with 
respect to $t\in [0,1]$ and $k\in {\Bbb N}$ by $h\in C^2(S)$. 
$I^{(2)}_{n_k}(t)$ is $P$-a.e. equicontinuous with respect to $t 
\in [0,1]$ and $k\in {\Bbb N}$ by the boundedness of $h'$ and 
Paul L\'evy's modulus of continuity for Brownian motion; modify, 
for example, the proof of \cite{MPSW10}, Theorem 1.12. In fact, we 
note that the modulus of continuity of $\int_{s=0}^\cdot\frac{1} 
{n_k}\sum_{i=1}^{n_k}h'\left(x^{n_k}_i(s)\right)\, d\beta_i(s)$ 
is majorized by the the modulus of continuity of $\|h'\|\frac{1} 
{n_k}\sum_{i=1}^{n_k}\beta_i$, $k\in {\Bbb N}$. The rest is just a 
slight modification of the calculation above (1.2) in \cite{MPSW10}. 

The term $I^{(3)}_{n_k}(t)$ is $P$-a.e. equicontinuous with respect  
to $t\in [0,1]$ and $k\in {\Bbb N}$. This follows from the fact that 
$\frac{1}{n_k}\sum_{i\neq j}\Psi\left(n_k^{1+\alpha}\left(x^{n_k 
}_i(t)-x^{n_k}_j(t)\right)\right)$ is $P$-a.e. uniformly bounded 
on $t\in [0,1]$ and $k\in {\Bbb N}$ according to Lemma 
\ref{Lemma4.1} (a) and $\Psi\le \beta V$, cf. introduction of this 
section. 

Let $h_1,h_2,\ldots\, \in C(S)$ be a sequence of linearly independent 
functions such that the collection of its finite linear combinations 
is dense in $C(S)$. 
Summarizing part (b) so far, we have shown, that $P$-a.e. there is 
a subsequence $n_{k_1}$ of $n_k$ such that $t\to (h_1,\mu_{n_{k_1}} 
(t,\cdot))-(h_1,\mu_{n_{k_1}}(0,\cdot))$ converges uniformly on 
$[0,1]$ to some $\rho(\cdot ,h_1)\in C[0,1]$ as $k_1\to\infty$. 
Iteratively, for $l\in {\Bbb N}$ there is a subsequence $n_{k_{l+1} 
}$ of $n_{k_l}$ such that $t\to (h_{l+1},\mu_{n_{k_{l+1}}}(t,\cdot) 
)-(h_{l+1},\mu_{n_{k_{l+1}}}(0,\cdot))$ converges uniformly on $[0, 
1]$ to some $\rho(\cdot ,h_{l+1})\in C[0,1]$ as $k_{l+1}\to\infty$. 
This holds $P$-a.e. simultaneously for all $l\in {\Bbb N}$ where 
the choice of the subsequences may depend on the element of 
$\Omega$. 

Thus, for $P$-a.e. elements of $\Omega$ there is a universal 
(diagonal) subsequence $n_q$, $q\in {\Bbb N}$, of $n_k$, $k\in {\Bbb 
N}$, such that for all $l\in{\Bbb N}$, $t\to (h_l,\mu_{n_q}(t,\cdot) 
)-(h_l,\mu_{n_q}(0,\cdot))$ converges uniformly on $[0,1]$ to $\rho( 
\cdot ,h_l)\in C[0,1]$ as $q\to\infty$.  
\medskip

\nid 
{\it Step 2 } We show the existence of a solution to (\ref{4.10}). 
For $P$-a.e. sequences of initial values $(x_1^{n_q}(0),\ldots , x_{ 
n_q}^{n_q}(0))$ and, respectively, initial empirical measures $\mu_{ 
n_q}(0,\cdot)$, $q\in {\Bbb N}$, we may choose a subsequence $n_r$, 
$r\in {\Bbb N}$, of $n_q$, $q\in {\Bbb N}$, and a measure $\mu(0, 
\cdot)$ on $(S,{\cal B}(S))$ such that $\mu_{n_r}(0,\cdot)\stack{r\to 
\infty}{\Ra}\mu(0,\cdot)$. 

For the result of Step 1 we replace the interval $t\in [0,1]$ by $t 
\in [0,T]$ for an arbitrary $T>0$. Since the linear hull of $h_1,h_2 
,\ldots\, $ in $C(S)$ is $C(S)$, for $P$-a.e. $\mu(0,\cdot)$ and 
every $t\in (0,T]$, 
\begin{eqnarray*}
\rho(t,h_l)+(h_l,\mu(0,\cdot))\, ,\quad l\in {\Bbb N}, 
\end{eqnarray*}
can be continuously extended to a linear functional $\mu^t(h)$ on 
$h\in C(S)$ such that 
\begin{eqnarray*}
\left(h,\mu_{n_r}(t,\cdot )\vphantom{l^1}\right)\stack{r\to\infty} 
{\lra}\mu^t(h)\quad\mbox{\rm on some subsequence $(n_r)_{r\in {\Bbb 
N}}$ of } (n_k)_{k\in {\Bbb N}} 
\end{eqnarray*}
for all $h\in C^2(S)$ where we mention once again that the choice of 
the subsequence $(n_r)_{r\in {\Bbb N}}$ may depend on the element of 
$\Omega$. Furthermore, $\mu^t(h)=\int_S h(\theta)\, \mu(t,d\theta)$, 
$h\in C(S)$, for some probability measure $\mu(t,\cdot )$. In 
particular we find 
\begin{eqnarray}\label{4.13}
\mu_{n_r}(t,\cdot)\stack{r\to\infty}{\Ra}\mu(t,\cdot)\, ,\quad t\ge 
0. 
\end{eqnarray}

Our task is now to demonstrate that $P$-a.e. 
\begin{eqnarray*}
\mu(t,d\theta)=\rho(t,\theta)\, d\theta\, ,\quad \theta\in S,\quad 
\mbox{\rm such that (\ref{4.7}) for $\mu(t,\cdot )$ and (\ref{4.10}) 
for all }t\ge 0. 
\end{eqnarray*}
Let $h\in C^2(S)$. We recall that on the subsequence $n_r$, $r\in 
{\Bbb N}$, (\ref{4.12}) reads as $\left(h,\mu_{n_r}(t,\cdot) 
\vphantom{l^1}\right)-\left(h,\mu_{n_r}(0,\cdot)\vphantom{l^1} 
\right)-I^{(1)}_{n_r}(t)-I^{(2)}_{n_r}(t)=I^{(3)}_{n_r}(t)$. By 
(\ref{4.13}) we have  
\begin{eqnarray}\label{4.14}
&&\hspace{-.5cm}\left(h,\mu_{n_r}(t,\cdot)\vphantom{l^1}\right)-
\left(h,\mu_{n_r}(0,\cdot)\vphantom{l^1}\right)-I^{(1)}_{n_r}(t) 
\nonumber \\ 
&&\hspace{.5cm}\stack{r\to\infty}{\lra}\left(h,\mu(t,\cdot) 
\vphantom{l^1}\right)-\left(h,\mu(0,\cdot)\vphantom{l^1}\right)- 
\frac12\int_{s=0}^t\left(h'',\mu(s,\cdot)\vphantom{l^1}\right)\, 
ds 
\end{eqnarray}
$P$-a.e. uniformly on $t\in [0,T]$ for all $T>0$. Next, let us 
examine $I^{(2)}_{n_r}$. Without loss of generality, we assume 
that $n_r\ge r^2$, $r\in {\Bbb N}$, and note that, for $N\in {\Bbb 
N}$, by 
\begin{eqnarray*}
&&\hspace{-.5cm}E\left(\sup_{r\in {\Bbb N},\, r\ge N}\left(I^{(2) 
}_{n_r}(T)\right)^2\right) \\ 
&&\hspace{.5cm}\le\sum_{r=N}^\infty\frac{1}{n_r^2} E\left(\int_{s 
=0}^t\sum_{i=1}^{n_r}h_1'\left(x^{n_r}_i(s)\right)\, d\beta_i(s) 
\right)^2\le\sum_{r=N}^\infty\frac{t\|h'\|}{n_r}\le\sum_{r=N 
}^\infty\frac{t\|h'\|}{r^2} \\ 
&&\hspace{.0cm}\stack{N\to\infty}{\lra}0\vphantom{\sum^\infty} 
\end{eqnarray*}
$\sup_{r\in{\Bbb N},\, r\ge N}(I^{(2)}_{n_r}(t))^2$, $t\ge 0$, is 
$P$-integrable and therefore a submartingale with respect to the 
filtration generated by $\beta_1,\beta_2,\ldots\, $. We verify now 
that 
\begin{eqnarray}\label{4.15}
I^{(2)}_{n_r}\, ,\ r\in{\Bbb N},\quad\mbox{\rm is $P$-a.e. uniformly 
bounded in }C([0,T])\ \mbox{\rm and }\lim_{r\to\infty}\left\|I^{(2) 
}_{n_r}\right\|_{C([0,T])}=0\qquad 
\end{eqnarray}
for every $T>0$ by using Doob's inequality and  
\begin{eqnarray*}
&&\hspace{-.5cm}P\left(\sup_{r\in {\Bbb N},\, r\ge N}\max_{t\in [0, 
T]}\left(I^{(2)}_{n_r}(t)\right)^2\ge a\right) \\ 
&&\hspace{.5cm}=P\left(\max_{t\in [0,T]}\sup_{r\in {\Bbb N},\, r\ge 
N}\left(I^{(2)}_{n_r}(t)\right)^2\ge a\right)\le\frac1aE\left(\sup_{ 
r\in {\Bbb N},\, r\ge N}\left(I^{(2)}_{n_r}(T)\right)^2\right)\, , 
\end{eqnarray*}
$a>0$. Let us turn to $I^{(3)}_{n_r}$. Lemma \ref{Lemma4.1} (a) 
together with $\Psi\le\beta V$ (cf. (\ref{4.3})) and Lemma 
\ref{Lemma4.2} together with (\ref{4.13}) imply 
\begin{eqnarray}\label{4.16}
I^{(3)}_{n_r}\stack{r\to\infty}{\lra}\int_{s=0}^tQ(\rho(s,\cdot),\cdot) 
(h')\, ds 
\end{eqnarray}
$P$-a.e. uniformly on $t\in [0,T]$ for all $T>0$. En passant we have 
also verified (\ref{4.7}) for all $t\ge 0$. 

The existence of a solution to (\ref{4.10}) is now a consequence of 
(\ref{4.12}) on the one hand and (\ref{4.14})-(\ref{4.16}) on the 
other hand. 
\medskip

\nid 
{\it Step 3 } We show uniqueness of the solution to (\ref{4.10}). Let 
us abbreviate 
\begin{eqnarray*}
L(h,\rho):=\frac12(\rho,h'')+Q(\rho,\cdot)(h') 
\end{eqnarray*}
for $\rho\in L^1(S)$ satisfying $\int_{\theta\in S}\int_{\tau\in S} 
|\theta-\tau|^{-\beta}|\rho(\theta)\rho(\tau)|\, d\tau\, d\theta\le 
B\|\rho\|^2_{L^1(S)}$ where $B$ is the constant from Lemma 
\ref{Lemma4.1}. 

Let $\rho_1(t,\theta)$, $t\ge 0$, and $\rho_2(t,\theta)$, $t\ge 0$, 
be two solutions to (\ref{4.10}) with $\rho_1(0,d\theta)=\rho_2(0, 
\theta)=\rho_0$ for some $\rho_0(\theta)\, d\theta$ satisfying 
(\ref{4.7}). Let $m_\ve(r,\cdot)$ be a usual family of one-dimensional 
mollifier functions, for fixed parameter $\ve>0$ symmetric about $r\in 
S$. Define 
\begin{eqnarray*}
\rho_{1,n}(\theta):=\rho_1\ast m_{\frac1n}(\theta)\equiv\int_{r\in 
S}\rho_1(r)m_{\frac1n}(r,\theta)\, dr\, ,\quad\theta\in S,\ n\in 
{\Bbb N}. 
\end{eqnarray*}
The same way define $\rho_{1,n}$. In addition, let $\rho^{2\ast}_{ 
1,n}:=\rho_{1,n}\ast m_{\frac1n}$, $\rho^{2\ast}_{2,n}:=\rho_{2,n} 
\ast m_{\frac1n}$. 

Assume the existence of $t_0>0$ such that $\left(\rho_1(t_0,\cdot) 
-\rho_2(t_0,\cdot),\rho_1(t,\cdot)-\rho_2(t,\cdot)\vphantom{l^1} 
\right)$ is increasing in $t_0$, i. e. the $\liminf$ of the 
differential quotient with respect to $t$ is positive in some 
neighborhood of $t_0$. Then there exists $n_0\in {\Bbb N}$ (large) 
and $t_1>0$ (near $t_0$) such that 
\begin{eqnarray*}
0&&\hspace{-.5cm}<\sup_{n\ge n_0}\left.\frac{d}{dt}\right|_{t=t_1} 
\left(\rho_{1,n}(t_1,\cdot)-\rho_{2,n}(t_1,\cdot),\rho_{1,n}(t,\cdot 
)-\rho_{2,n}(t,\cdot)\vphantom{l^1}\right)=:\delta\, .\vphantom 
{\int_0^t} 
\end{eqnarray*}
For the existence of the derivative recall (\ref{4.10}) and that 
$\rho_{1,n}(t,\theta)$ has the form $(\rho_1(t,\cdot),h)$ with $h= 
m_{\frac1n}(\cdot,\theta)$. The same holds for $\rho_{2,n}(t,\theta 
)$. We choose $n\ge n_0$ and $\alpha,\beta>0$ such that with $\rho:= 
\alpha\left(\rho_{1,n}(t_1,\cdot)-\rho_{2,n}(t_1,\cdot)\right)+\beta$ 
we have 
\begin{eqnarray*} 
\t \rho:=\alpha\left(\rho_{1,n}^{2\ast}(t_1,\cdot)-\rho_{2,n}^{2\ast} 
(t_1,\cdot)\right)+\beta\ge 0\, ,\quad\int_{\theta\in S}\t\rho(\theta 
)\, d\theta=1\, ,  
\end{eqnarray*}
and 
\begin{eqnarray}\label{4.17}
\left|\frac{1}{\alpha}\left(L\left(\t \rho,\rho_1(t_1,\cdot)\right) 
-L\left(\t \rho,\rho_2(t_1,\cdot)\right)\vphantom{l^1}\right)-\frac 
{1}{\alpha^2}L(\t \rho,\rho)\right|<\frac{\delta}{4}\, ,  
\end{eqnarray}
as well as 
\begin{eqnarray}\label{4.18}
\left|\frac{1}{\alpha^2}L(\t \rho,\rho)-\frac{1}{\alpha^2}L(\t \rho, 
\t \rho)\right|<\frac{\delta}{4}\, ; 
\end{eqnarray}
note that, for $\beta\to\infty$ and $\alpha$ accordingly adjusted, 
the left-hand side of (\ref{4.17}) tends to zero and that 
(\ref{4.18}) can be achieved by choosing $n$ sufficiently large. 
We have 
\begin{eqnarray}\label{4.19}
0&&\hspace{-.5cm}<\delta=\frac{1}{\alpha^2}\left.\frac{d}{dt}\right 
|_{t=t_1}\left(\alpha\left(\rho_{1,n}(t_1,\cdot)-\rho_{2,n}(t_1, 
\cdot)\right)+\beta,\alpha\left(\rho_{1,n}(t,\cdot)-\rho_{2,n}(t, 
\cdot)\right)+\beta\vphantom{l^1}\right)\nonumber \\ 
&&\hspace{-.5cm}=\frac{1}{\alpha}L\left(\t \rho,\rho_1(t_1,\cdot) 
\right)-\frac{1}{\alpha}L\left(\t \rho,\rho_2(t_1,\cdot)\right) 
<\frac{1}{\alpha^2}L(\t \rho,\t \rho)+\frac{\delta}{2}\, . 
\end{eqnarray}

Next we aim to show that the right-hand side of (\ref{4.19}) does 
not exceed $\delta/2$ which will show the above assumption does 
not hold. We have thus proved uniqueness. Let 
\begin{eqnarray*} 
\t \mu_{n'}\equiv\frac{1}{n'}\sum_{i=1}^{n'}\delta_{\t x_{i,{n'}}} 
\in {\cal M}_1^{n'}(S)\quad\mbox{\rm such that}\quad\left(\t \rho,\t 
\mu_{n'}\right)=\max_{\mu_{n'}\in {\cal M}_1^{n'}(S)}(\t \rho,\mu_{n 
'})\, . 
\end{eqnarray*}
By the maximum principle for infinitesimal operators we have 
\begin{eqnarray}\label{4.20} 
{\cal L}_{n'}\left(\t \rho,\t \mu_{n'}\right)\equiv{\cal L}_{n'}\left 
(\t \rho,\frac{1}{n'}\sum_{i=1}^{n'}\delta_{\t x_{i,{n'}}}\right)\le 
0\, ,\quad n'\in {\Bbb N}. 
\end{eqnarray}
Furthermore, $\t \mu_{n'}\stack{{n'}\to\infty}{\Ra}\t \rho(\theta) 
\, d\theta$ and therefore 
\begin{eqnarray*} 
{\cal L}_{n'}\left(\t \rho,\t \mu_{n'}\right)\equiv{\cal L}_{n'} 
\left(\t \rho,\frac{1}{n'}\sum_{i=1}^{n'}\delta_{\t x_i}\right) 
\stack{{n'}\to\infty}{\lra}L(\t\rho,\t\rho)\, . 
\end{eqnarray*}
Together with (\ref{4.19}) and (\ref{4.20}) this completes the proof 
of uniqueness. 
\medskip 

\nid 
(c) This is just a modification of Step 3 (uniqueness) of part (b). 
\medskip 

\nid 
(d) Let $\hat{\rho}(\theta)\, d\theta\neq\mu_0(d\theta)$ be a 
probability measure with (\ref{4.7}) and let $\hat{\rho}(t,\cdot)$, 
$t\ge 0$, be the solution to (\ref{4.10}) with $\hat{\rho}(0,\cdot):= 
\hat{\rho}$. There is a version $\hat{\rho}$ and a maximum point 
$\theta_0\in S$ of this version $\hat{\rho}$ in the sense that ess$ 
\, \sup_{\theta\in U}\hat{\rho}(\theta)=\hat{\rho}(\theta_0)$ for all 
open $U\subset S$ with $\theta_0\in U$. We note $\hat{\rho}(\theta_0) 
>1$. Without loss of generality, we may assume that $\theta_0=0$. Let 
$\check{\rho}$ be defined by $\check{\rho}(\theta):=\hat{\rho}(-\theta 
)$, $\theta\in S$. Furthermore, let $\check{\rho}(t,\cdot)$, $t\ge 0$, 
be the solution to (\ref{4.10}) with $\check{\rho}(0,\cdot):=\check{ 
\rho}$. Considering only $h\in C^2(S)$ which are symmetric about 
$\theta_0=0$ we observe 
\begin{eqnarray}\label{4.21}
\int h''(\theta)\hat{\rho}(\theta)\, d\theta=\int h''(\theta)\check 
{\rho}(\theta)\, d\theta\quad\mbox{\rm and}\quad Q\left(\hat{\rho}, 
\cdot\right)(h')=-Q\left(\check{\rho},\cdot\right)(h')\, . 
\end{eqnarray} 
On the other hand, by symmetry of the system, 
\begin{eqnarray}\label{4.22} 
\quad\frac{d}{dt}\left(\t h,\hat{\rho}(t,\cdot)\right)=0\quad\mbox 
{\rm for all $t\ge 0$ and all $\t h\in C^2(S)$ if and only if}\quad 
\frac{d}{dt}\left(\t h,\check{\rho}(t,\cdot)\right)=0\nonumber \\ 
\end{eqnarray} 
for all $t\ge 0$ and all $\t h\in C^2(S)$. Assuming (\ref{4.11}) for 
$\hat{\rho}$ then we have (\ref{4.22}) and thus (\ref{4.11}) also for 
$\check{\rho}$. In particular, we get 
\begin{eqnarray*}  
\left.\frac{d}{dt}\right|_{t=0}\left(h,\hat{\rho}(t,\cdot)\right)= 
0\quad\mbox{\rm and}\quad\left.\frac{d}{dt}\right|_{t=0}\left(h, 
\check{\rho}(t,\cdot)\right)=0 
\end{eqnarray*} 
and therefore 
\begin{eqnarray*}  
\frac12\int h''(\theta)\hat{\rho}(\theta)\, d\theta+Q\left(\hat{\rho} 
,\cdot\right)(h')=0\quad\mbox{\rm and}\quad\frac12\int h''(\theta) 
\check{\rho}(\theta)\, d\theta+Q\left(\check{\rho},\cdot\right)(h')=0
\end{eqnarray*} 
for all $h\in C^2(S)$ which are symmetric about $\theta_0=0$ which is 
a maximum point in the above sense of $\hat{\rho}$ as well as $\check 
{\rho}$. This yields a contradiction to (\ref{4.21}). 
\qed
\medskip 

Now we turn to initial measures $\bnu_n$ for $(X^n_t)_{t\ge 0}$, $n\in 
{\Bbb N}$, the invariant measures. We recall that $\mu_0$ is the measure 
on $(D,{\cal B}(D))$ which is the uniform distribution on $D=S$. 

Let us now introduce the stationary version ${\bf X}^n=\left((X^n_t)_{t 
\ge 0},P^n\right)$ defined by $P^n:=\int P^n_\mu\, {\bnu}_n(d\mu)$, $n 
\in {\Bbb N}$. Let ${\bf X}$ be the path concentrated on $\mu_0$.  As in 
Section 3, let $(T_{n,t})_{t\ge 0}$ denote the semigroup associated with 
$X^n$. The measure ${\bnu}_n$ is an invariant measure of the semigroup 
$(T_{n,t})_{t\ge 0}$, $n\in {\Bbb N}$. This follows from the definition 
of the measures $P^n_\mu$, $\mu\in {\cal M}_1^n(D)$, and the fact that 
the measure $\nu_n$ is invariant for the diffusion $x^n=((x^n(t))_{t\ge 
0},\hat{P}^n_\eta)$,   
\begin{eqnarray*} 
\int P^n_\mu(X^n_t\in {\cal A}(A))\, {\bnu}_n(d\mu )&&\hspace{-.5cm}= 
\int\hat{P}^n_\eta(x^n_t\in A)\, \nu_n(d\eta ) \\ 
&&\hspace{-.5cm}=\nu_n (A)={\bnu}_n({\cal A}(A))\, , \quad A\in {\cal B 
}^{\Pi (n)}(S). \vphantom{\sum}
\end{eqnarray*} 

Let ${\bnu}$ be the probability measure on $(E,{\cal B}(E)):=({\cal M 
}_1(D),{\cal B}({\cal M}_1(D)))$ which is concentrated on $\mu_0$. Our 
goal is to prove the following theorem. 
\begin{theorem}\label{Theorem4.4}  
(a) We have ${\bnu}_n\stack{n\to\infty}{\Ra}{\bnu}=\delta_{\mu_0}$. \\ 
(b) The processes ${\bf X}^n$ converge weakly to ${\bf X}$ as $n\to 
\infty$. 
\end{theorem}

We choose a sequence of linearly independent functions $h_1,h_2, 
\ldots\, \in C^\infty(S)$ such that the collection of its finite 
linear combinations is dense in $C(S)$. 
Moreover we will work with the space $\t C_b(E)$ specified by 
\begin{eqnarray}\label{4.23}
&&\hspace{-1.6cm}\t C_b(E):=\left\{g(\mu)=\vp((h_1,\mu),(h_2,\mu), 
\ldots\, ),\ \mu\in E:\vphantom{\textstyle\sum_{i,j=1}^\infty\left|\
frac{\partial^2\vp}{\partial y_i\partial y_j}\right|}\vp\in C_b^2 
({\Bbb R}^{\Bbb N}),\ \vphantom{\sum_{i,j=1}^\infty}\right.\nonumber 
 \\ 
&&\hspace{-1.2cm}\left. 
\sum_{i=1}^\infty\left\|\frac{\partial\vp}{\partial y_i}\right\|\left( 
1+\|\nabla h_i\|+\|\Delta h_i\|\right)\le c,\ \sum_{i,j=1}^\infty\left 
\|\frac{\partial^2\vp}{\partial y_i\partial y_j}\right\|\left(\|\nabla 
h_i\|\|\nabla h_j\|\right)\le c,\ c>0\right\}. \nonumber \\  
\end{eqnarray}
Referring to \cite{Lo05}, Subsection 2.7, define 
\begin{eqnarray}\label{4.24}
&&\hspace{-.5cm}\frac12\Del g(\mu):=\frac12\sum_{i=1}^\infty\frac{ 
\partial\vp (\mu)}{\partial x_i}\cdot (\Delta h_i,\mu)+\frac{1}{2n} 
\sum_{i,j=1}^\infty\frac{\partial^2\vp (\mu)}{\partial x_i\partial x_j} 
\cdot (\nabla h_i\circ \nabla h_j,\mu)\, , \nonumber\\ 
&&\hspace{8cm}\mu\in E_n\, , \ n\ge 2\, , \ g\in \t C_b(E),  
\end{eqnarray}
$\circ$ indicates the scalar product in ${\Bbb R}^d$, and $\nn\frac12 
\Del g\nn :=\sup_{\mu\in \bigcup_{n\ge 2}E_n}\|\frac12\Del g(\mu)\|$. 
Furthermore, we have used $\frac{\partial\vp(\mu)}{\partial y_k}$ as an 
abbreviation for $\frac{\partial\vp}{\partial y_k}((h_1,\mu),\ldots \, 
)$ and $\frac{\partial^2\vp (\mu)}{\partial x_i\partial x_j}$ as an 
abbreviation for $\frac{\partial^2\vp}{\partial x_i\partial x_j}((h_1, 
\mu),\ldots \, )$. From (\ref{4.5}), we obtain for $g\in\t C_b(E)$ and 
$\mu =\frac1n\sum_{i=1}^n\delta_{x_i}$ 
\begin{eqnarray}\label{4.25}
A_n g(\mu)=\frac12\Del g(\mu)-\sum_{l=1}^\infty\frac{\partial\vp(\mu)} 
{\partial y_l}\sum_{j\neq i}n^\alpha V'(n^{1+\alpha}(x_i-x_j))h_l'(x_i) 
\, , \quad n\in{\Bbb N}. 
\end{eqnarray}
Note that the factor $n^{1+\alpha}$ in the second item of (\ref{4.5}) 
reduces here to $n^\alpha$ by the factor $\frac1n$ coming from $\frac 
{\partial}{\partial x_i}\sum_{j=1}^n\frac1n h_l(x_j)$. Let us also 
define $(T_t)_{t\ge 0}$ regarded as a semigroup in $L^2(E,\t\bnu)$ 
where $\t\bnu$ is a accumulation point of $({\bnu}_n)_{n\in {\Bbb N}}$ 
by 
\begin{eqnarray*}
T_tf(\rho(\theta)\, d\theta):=\rho(t,\theta)\, d\theta\, ,\quad 
\rho(\theta)\, d\theta\equiv\rho(0,\theta)\, d\theta\quad\mbox{\rm 
satisfies (\ref{4.7})} 
\end{eqnarray*}
where we recall Lemma \ref{Lemma4.1} (b) and Proposition 
\ref{Proposition4.3} (b). That $T_tf\in L^2(E,\t\bnu)$ if $f\in L^2(E, 
\t\bnu)$, $t\ge 0$, follows from the fact that $\t\bnu$ is an invariant 
measure, cf. Remark (9) of Section 2. 
\medskip

In the following proposition, we collect all the necessary prerequisites 
in order to deduce Theorem \ref{Theorem4.4} (a) from Proposition 
\ref{Proposition3.4} (a) and Theorem \ref{Theorem4.4} (b) from 
Corollary \ref{Corollary3.6} together with Theorem \ref{Theorem4.4} (a). 
\begin{proposition}\label{Proposition4.5} 
(a) $(T_{n,t})_{t\ge 0}$ is a strongly continuous semigroup in the space 
$L^2(E_n,{\bnu}_n)$. \\ 
(b) The conditions of Proposition \ref{Proposition3.4} (a) are satisfied, 
i. e., we have Theorem \ref{Theorem4.4} (a). \\ 
(c) We have (${\cal C}4$) for $\beta\ge 1$. 
\end{proposition}
Proof. (a) Given $F\in L^2(E_n,{\bnu}_n)$, let $f\in L^2(S^n,\nu_n)$ be 
the function symmetric in all variables satisfying $F=\t F_{f,n}$. From 
the definitions of the measures $P^n_\mu$, $\mu\in {\cal M}_1^n(D)$, and 
the semigroup $(T_{n,t})_{t\ge 0}$, we obtain 
\begin{eqnarray*} 
\|T_{n,t}F-F\|_{L^2(E_n,\sbnu_n)}&=&\left\|\int F(\nu)\, P^n_\cdot 
(X_t\in d\nu)-F\right\|_{L^2(E_n,\sbnu_n)} \\ 
&=&\left\|\int f(y)\, \hat{P}^n_\cdot(x_t\in dy)-f\right\|_{L^2(S^n, 
\nu_n)}=\|\hat{T}_{n,t}f-f\|_{L^2(S^n,\nu_n)}\, . 
\end{eqnarray*}
The claim is now a consequence of the fact that $(\hat{T}_{n,t})_{t\ge 
0}$ forms a strongly continuous semigroup on $L^2( S^n,\nu_n)$. 
\medskip 

\nid 
(b) {\it Step 1 } We will use Proposition \ref{Proposition3.4} (a). Let 
us assume that we have fixed $\t\bnu$ and that we have chosen $\hat{C} 
(E):=C_b(E)$. At the same time let us review the Remarks (1) and (2) of 
Section 3. 

In this step, let us verify the conditions (i), (ii), and (v) of 
Proposition \ref{Proposition3.4} and show that $(T_t)_{t\ge 0}$ regarded 
as a semigroup in $L^2(E,\t\bnu)$ is strongly continuous. 

Condition (${\cal C}3$) is satisfied by the definition given immediately 
before Theorem \ref{Theorem4.4} according to which $\bnu=\delta_{\mu_0}$. 
By using Proposition \ref{Proposition4.3} (d) we verify (v) of 
Proposition \ref{Proposition3.4}. For condition (i) of Proposition 
\ref{Proposition3.4} we note that $(T_{n,t})_{t\ge 0}$ is Feller; for 
this see also E. B. Dynkin \cite{Dy65}, Theorem 5.11. This shows $\{G_{n, 
\beta}g: g\in \hat{C}_b(E),\ \beta >0\}\subseteq\hat{C}_b(E)$ in the 
sense of (i) of Proposition \ref{Proposition3.4}. For condition (ii) of 
Proposition \ref{Proposition3.4} we refer to Proposition 
\ref{Proposition4.3} (c) to show that $(T_t)_{t\ge 0}$ is Feller. 

By Lemma \ref{Lemma4.1} (b) and Proposition \ref{Proposition4.3} (b), 
we have for $\t\bnu$-a.e. $\mu\equiv\rho(\theta)\, d\theta\equiv\rho (0, 
\theta)\, d\theta$ 
\begin{eqnarray*} 
&&\hspace{-.5cm}\left|\int J(\theta )\rho(t,\theta)\, d\theta-\int J 
(\theta)\rho(\theta)\, d\theta\right|=\left|\int_0^t\left(\frac12\int 
J''(\theta)\rho(s,\theta)\, d\theta+Q(\rho(s,\cdot),\cdot)(J')\right) 
\, ds\right| \\ 
&&\hspace{.5cm}\le\left\|J''\right\|\left(\frac{1}{2}+B\right)t\, , 
\quad J\in C^\infty(S),  
\end{eqnarray*}
where $B$ is the constant from Lemma \ref{Lemma4.1} (a). From here we 
get 
\begin{eqnarray}\label{4.26} 
\int |T_tg-g|\, d\t\bnu\stack{t\to 0}{\lra}0\, ,\quad g\in C_b(E),  
\end{eqnarray}
and thus $\left(\int (T_tg-g)^2\, d\t\bnu\right)^{1/2}\le\left(2\|g\| 
\right)^{1/2}\left(\int |T_tg-g|\, d\t\bnu\right)^{1/2}\stack{t\to 0} 
{\lra}0$. In order to show that $(T_t)_{t\ge 0}$ regarded as a semigroup 
in $L^2(E,\t\bnu)$ is strongly continuous, we note first that by 
(\ref{4.26}) and Proposition \ref{Proposition4.3} (c) $(T_t)_{t\ge 0}$ 
regarded as a semigroup in $C_b(E)$ is strongly continuous. Denote this 
semigroup in $C_b(E)$ by $(T^c_t)_{t\ge 0}$ and its generator by $A^c$. 
We know that the domain $D(A^c)$ is dense in $L^2(E,\t\bnu)$. From 
$T^c_tg-g=\int_0^tT^c_sA^cg\, ds$, we conclude 
\begin{eqnarray}\label{4.27}
\|T^c_tg-g\|\le t\cdot \|A^cg\|\, , \quad g\in D(A^c). 
\end{eqnarray}
Now, let $f\in L^2(E,\t\bnu)$, $\ve >0$, and choose $g\in D(A^c)$ such 
that $\|f-g\|_{L^2(E,\t\sbnu)}<\frac{\ve}{3}$. Furthermore, choose $t< 
\ve/(3\|A^cg\|)$. Because of (\ref{4.27}) we then have $\|T^c_tg-g\|< 
\frac{\ve}{3}$. In addition, it holds that $\|T_t(f-g)\|_{L^2(E,\t\sbnu)}
\le\|f-g\|_{L^2(E,\t\sbnu)}<\frac{\ve}{3}$ where we have employed the 
contractivity of the semigroup $(T_t)_{t\ge 0}$ in $L^2(E,\t \bnu)$ 
as a consequence of the fact that $\t\bnu$ is an invariant measure for 
$(T_t)_{t\ge 0}$, cf. also Remark (9) of Section 2. Finally, we obtain 
\begin{eqnarray*}
\left\|T_tf-f\right\|_{L^2(E,\t\sbnu)}&&\hspace{-.5cm}\le\|T_t(f-g)\|_{ 
L^2(E,\t\sbnu)}+\|T_tg-g\|_{L^2(E,\t\sbnu)}+\|g-f\|_{L^2(E,\t\sbnu)} \\ 
&&\hspace{-.5cm}\le\textstyle\frac{\ve}{3}+\|T^c_tg-g\|+\textstyle\frac 
{\ve}{3}<\ve\, . 
\end{eqnarray*}
Thus, $(T_t)_{t\ge 0}$ is strongly continuous in $L^2(E,\t\bnu)$. 
\medskip 

\nid 
{\it Step 2 } The next step is devoted to the verification of 
condition (iii) of Proposition \ref{Proposition3.4}. The verification of 
(iv) of Proposition \ref{Proposition3.4} will be en passant. 

For an arbitrary accumulation point $\t\bnu$ of $(\bnu_n)_{n\in {\Bbb 
N}}$ and a subsequence $n_k$ such that $\bnu_{n_k}\stack{k\to\infty}{\Ra} 
\t\bnu$ we will show that $S_{n_k}$ converges to $S_{\t \sbnu}$ in the 
sense of Remark (7) of Section 2 with ${\cal C}$ as in Definition 
\ref{Definition3.1}. We are going to verify (ii') and (iii'') of Remark 
(7) of Section 2. Let us first assume that $\psi\equiv g\in\t C_b(E)$ as 
given in (\ref{4.23}). It is standard that $\psi\in D(A_n)$, $n\in {\Bbb 
N}$. Thus, we can choose $\Gamma:=\t C_b(E)$ in (iv) of Proposition 
\ref{Proposition3.4}. Let $\t\psi\in\t C_b(E)$. 

Using the abbreviations $\frac{\partial\vp}{\partial y_k}$ for $\frac{ 
\partial\vp}{\partial y_k}(\frac1n\sum_{i=1}^nh_1(x_i),\ldots\, )$ and 
$\t\psi$ for $\t\psi(\frac1n\sum_{i=1}^n\delta_{x_i})$, we have 
\begin{eqnarray}\label{4.28}
&&\hspace{-.5cm}S_{n_k}(\psi,\t\psi)=-\langle A_{n_k}\psi,\t\psi 
\rangle_{n_k}\nonumber \\ 
&&\hspace{.0cm}=-\frac12\langle\Del_{n_k}\psi,\t\psi\rangle_{n_k}+ 
n_k^\alpha\sum_{l=1}^\infty\sum_{j\neq i}\int_{S^{n_k}}\t\psi\frac 
{\partial\vp}{\partial y_l}\cdot V'(n_k^{1+\alpha}(x_i-x_j))h_l' 
(x_i)\, \nu_{n_k}(dx)\nonumber \\ 
&&\hspace{.0cm}=-\frac12\langle\Del_{n_k}\psi,\t\psi\rangle_{n_k}+ 
\frac{n_k^\alpha}{2}\sum_{l=1}^\infty\sum_{j\neq i}\int_{S^{n_k}} 
\t\psi\frac{\partial\vp}{\partial y_l}\cdot\left(h_l'(x_i)-h_l'(x_j) 
\right)V'(n_k^{1+\alpha}(x_i-x_j))\, \nu_{n_k}(dx)\nonumber \\ 
&&\hspace{.0cm}=-\frac12\langle\Del_{n_k}\psi,\t\psi\rangle_{n_k}- 
\frac{1}{2n_k}\sum_{l=1}^\infty\sum_{j\neq i}\int_{S^{n_k}}\t\psi 
\frac{\partial\vp}{\partial y_l}\cdot\frac{h_l'(x_i)-h_l'(x_j)}{x_i 
-x_j}\Psi(n_k^{1+\alpha}(x_i-x_j))\, \nu_{n_k}(dx)\, .\nonumber \\ 
\end{eqnarray}
Also by (\ref{4.24}) and (\ref{4.23}), 
\begin{eqnarray}\label{4.29} 
-\frac12\langle\Del_{n_k}\psi,\t\psi\rangle_{n_k}\stack{k\to\infty} 
{\lra}-\frac12\int\sum_{l=1}^\infty\frac{\partial\vp (\mu)}{\partial 
y_l}\cdot (\Delta h_l,\mu)\t\psi (\mu)\, \t \bnu (d\mu)\, . 
\end{eqnarray}
Using the abbreviations $\left(\t\psi\frac{\partial\vp}{\partial y_l}
\right)(x)$ for $\t\psi(\frac1n\sum_{i=1}^n\delta_{x_i})\cdot\frac{ 
\partial\vp}{\partial y_l}(\frac1n\sum_{i=1}^nh_1(x_i),\ldots\, )$ and 
$\left(\t\psi\frac{\partial\vp}{\partial y_l}\right)\left(\rho\right)$ 
for $\t\psi(\rho(\theta)\, d\theta)\cdot\frac{\partial\vp}{\partial 
y_l}\left(\int h_1(\theta)\rho(\theta)\, d\theta,\ldots\, \right)$ 
where $\mu=\rho(\theta)\, d\theta$ it follows that 
\begin{eqnarray}\label{4.30}
&&\hspace{-.5cm}-\frac{1}{2n_k}\sum_{l=1}^\infty\sum_{j\neq i}\int_{ 
S^{n_k}}\left(\t\psi\frac{\partial\vp}{\partial y_l}\right)(x)\cdot 
\frac{h_l'(x_i)-h_l'(x_j)}{x_i-x_j}\Psi(n_k^{1+\alpha}(x_i-x_j))\, 
\nu_{n_k}(dx) \nonumber \\ 
&&\hspace{.0cm}\stack{k\to\infty}{\lra}\hspace{-.4cm}-\frac{\beta}{2} 
\sum_{l=1}^\infty\int\left(\t\psi\frac{\partial\vp}{\partial y_l} 
\right)\left(\rho\right)\cdot\int_{\theta\in S}\int_{\tau\in S}\frac 
{h_l'(\theta)-h_l'(\tau)}{\theta-\tau}|\theta-\tau|^{-\beta}\rho 
(\theta)\rho(\tau)\, d\tau\, d\theta\, \t \bnu(d\mu)\nonumber 
 \\ 
&&\hspace{.5cm}=-\sum_{l=1}^\infty\int\left(\t\psi\frac{\partial\vp} 
{\partial y_l}\right)\left(\rho\right)\cdot Q(\rho,\cdot)(h_l')\, \t 
\bnu(d\mu)
\end{eqnarray}
where, for the convergence, we have applied Lemma \ref{Lemma4.2} (b). 
Using now (\ref{4.10}), the notation $\rho(0,\theta)=\rho(\theta)$, 
and the right derivative 
\begin{eqnarray}\label{4.31}
&&\hspace{-.5cm}\left.\frac{d}{dt}\right|_{t=0}\psi(\rho(t,\theta)\, 
d\theta)\nonumber \\ 
&&\hspace{.5cm}=\sum_{l=1}^\infty\frac{\partial\vp (\rho(\theta) 
\, d\theta)}{\partial y_l}\cdot\left(\frac12\int_{\theta\in S}h_l''( 
\theta)\rho(\theta)\, d\theta+Q(\rho,\cdot)(h_l')\right)\quad\t \bnu 
\mbox{\rm -a.e.},  
\end{eqnarray}
we get from (\ref{4.30})
\begin{eqnarray}\label{4.32} 
&&\hspace{-.5cm}-\frac{1}{2n_k}\sum_{l=1}^\infty\sum_{j\neq i}\int_{ 
S^{n_k}}\left(\t\psi\frac{\partial\vp}{\partial y_l}\right)(x)\cdot 
\frac{h_l'(x_i)-h_l'(x_j)}{x_i-x_j}\Psi(n_k^{1+\alpha}(x_i-x_j))\, 
\nu_{n_k}(dx)\nonumber \\ 
&&\hspace{.0cm}\stack{k\to\infty}{\lra}-\int\left.\frac{d}{dt} 
\right|_{t=0}\psi(\rho(t,\theta)\, d\theta)\cdot\t\psi (\mu)\, \t \bnu 
(d\mu)+\frac12\int\sum_{l=1}^\infty\frac{\partial\vp (\mu)}{\partial 
y_l}\cdot (\Delta h_l,\mu)\t\psi (\mu)\, \t\bnu (d\mu)\, .\nonumber \\   
\end{eqnarray}
As in (\ref{4.24}) we have used the notation $\frac{\partial\vp (\mu)} 
{\partial y_l}$ for $\frac{\partial\vp}{\partial y_l}((h_1,\mu),\ldots\, 
)$. It follows from Lemma \ref{Lemma4.1} (b) 
that  
\begin{eqnarray}\label{4.33} 
\frac1t\int\int_{s=0}^t\left(\frac12\int_{\theta\in S}h''(\theta)\rho 
(s,\theta)\, d\theta+Q(\rho(s,\cdot),\cdot)(h')\right)^2\, ds\, d\t 
\bnu\le\|h''\|^2\left(\frac{1}{2}+B\right)^2\, , 
\end{eqnarray}
$h\in C^2(S)$, where $B$ is the constant of Lemma \ref{Lemma4.1}. 
According to representation (\ref{4.31}), (\ref{4.33}), and the 
definition of $\t C_b(E)$ the right derivative (\ref{4.31}) exists in 
$L^2(E,\t \bnu)$ for $\psi\in\t C_b(E)$ and we have 
\begin{eqnarray*} 
A_{\t\sbnu}\psi(\rho(\theta)\, d\theta)=\left.\frac{d}{dt}\right|_{t=0} 
\psi (\rho(t,\theta)\, d\theta)\quad\mbox{\rm in } L^2(E,\t \bnu)\, .  
\end{eqnarray*}
Furthermore, $\t C_b(E)\subseteq D(A_{\t\sbnu})$. Thus, from (\ref{4.32}) 
we obtain 
\begin{eqnarray}\label{4.34} 
&&\hspace{-.5cm}-\frac{1}{2n_k}\sum_{l=1}^\infty\sum_{j\neq i}\int_{ 
S^{n_k}}\left(\t\psi\frac{\partial\vp}{\partial y_l}\right)(x)\cdot 
\frac{h_l'(x_i)-h_l'(x_j)}{x_i-x_j}\Psi(n_k^{1+\alpha}(x_i-x_j))\, 
\nu_{n_k}(dx)\nonumber \\ 
&&\hspace{.5cm}\stack{k\to\infty}{\lra}-\frac12\int A_{\t \sbnu}\psi 
\cdot\t\psi\, d\t \bnu +\frac12\int\sum_{l=1}^\infty\frac{\partial\vp 
(\mu)}{\partial y_l}\cdot (\Delta h_l,\mu)\t\psi(\mu)\, \t \bnu (d\mu) 
\, . 
\end{eqnarray}
Piecing (\ref{4.28})-(\ref{4.32}), (\ref{4.34}) together we arrive at 
\begin{eqnarray}\label{4.35} 
S_{n_k}(\psi,\t\psi)\stack{k\to\infty}{\lra}S_{\t \sbnu}(\psi,\t\psi)\, , 
\quad\psi,\t\psi\in\t C_b(E). 
\end{eqnarray}
Calculations similar to (\ref{4.28})-(\ref{4.32}) also show that 
\begin{eqnarray}\label{4.36} 
\langle A_{n_k}\psi,A_{n_k}\psi\rangle_{n_k}\stack{k\to\infty}{\lra} 
\|A_{\t \sbnu}\psi\|_{L^2(E,\t \sbnu)}^2\, , \quad\psi\in\t C_b(E). 
\end{eqnarray}
For $\hat{\mu}\in E$, the function 
\begin{eqnarray*}
H_{\hat{\mu}}:=\sum_{k=1}^\infty {\frac{2^{-k}}{1+(\|h_k\|+\|\nabla 
h_k\|)^2+\|h_k\Delta h_k\|}}\left((h_k,\cdot)-(h_k,\hat{\mu})\right)^2
\end{eqnarray*}
belongs to $\t C_b(E)$ and $\{H_{\hat{\mu}}:\hat{\mu}\in E\}$ 
separates the points in $E$. Since $\t C_b(E)$ given by (\ref{4.23}) 
forms an algebra containing the constant functions the set $\t C_b(E)$ 
is dense in $C_b(E)$ and thus in $L^2(E,\t\bnu)$. 
\newpage

By (\ref{4.35}) $(S_{\t\sbnu},\t C_b(E))$ is symmetric and positive. 
By the Friedrichs extension, $(S_{\t\sbnu},\t C_b(E))$ is closable 
and its closure has a self-adjoint generator. This generator is  
$(A_{\t\sbnu},D(A_{\t \sbnu}))$ as a consequence of, for example, 
(\ref{4.35}) and \cite{P83}, Chapter 1, Corollary 4.4. The latter 
reference says that a densely defined closed linear operator, which 
together with its adjoint is dissipative, is the generator of a 
strongly continuous semigroup. 
\medskip

By the self-adjointness of $A_{\t\sbnu}$ and (\ref{4.34}) it follows 
that 
\begin{eqnarray*}
\int A_{\t\sbnu}\psi\cdot\t\psi\, d\t\bnu=\frac12\int A_{\t\sbnu} 
(\psi\t\psi)\, d\t\bnu=\frac12\int (\psi\t\psi)\cdot A_{\t\sbnu}\1 
\, d\t\bnu=0\, , \quad \psi,\t\psi\in \t C_b(E), 
\end{eqnarray*}
i.e., $A_{\t\sbnu}\psi =0$, $\psi\in\t C_b(E)$. With (\ref{4.35}) and 
(\ref{4.36}) we get the convergence to zero in (iv) of Proposition 
\ref{Proposition3.4} (a). The associated semigroup in $L^2(E,\t\bnu)$ 
is $T_t=identity$, $t\ge 0$. Thus, $D(A_{\t\sbnu})=D(S_{\t\sbnu})=L^2 
(E,\t\bnu)$. 

For $\psi\in D(S)$, there is a subsequence $n_r$, $r\in {\Bbb N}$, 
of $n_k$, $k\in {\Bbb N}$, and a sequence $\psi_{n_r}\in\t C_b(E)$, 
$r\in {\Bbb N}$, such that $\psi_{n_r}$ $s$-converges to $\psi$ and 
converges to $\psi$ in $L^2(E,\t\bnu)$ as $r\to\infty$, cf. 
Proposition \ref{Proposition3.3} (c). Furthermore, by a slight 
modification of Proposition \ref{Proposition3.3} (c) motivated by 
(\ref{4.35}) (add in the third sentence of its proof the line 
$-\langle A_m\t \vp_n,\t \vp_n\rangle_m\stack{m\to\infty}{\lra}0$ 
and proceed accordingly) it follows that 
\begin{eqnarray*} 
S_{n_r}(\psi_{n_r},\psi_{n_r})\stack{r\to\infty}{\lra}S_{\t 
\sbnu}(\psi_{n_r},\psi_{n_r})=0\, . 
\end{eqnarray*}
We get (ii') of Remark (7) of Section 2. Condition (i) of Definition 
\ref{Definition2.4} is now trivial. The bilinear forms $S_n$, $n\in 
{\Bbb N}$, and $S$ are symmetric in the sense of Remark (4) of Section 
2. This yields (iii'') of Remark (7) of Section 2. Summing up, we 
have verified (iii) of Proposition \ref{Proposition3.4} (a). 
\medskip 

\nid 
(c) Observe as above that $\t C_b(E)$ given by (\ref{4.23}) forms an 
algebra containing the constant functions and separating the points in 
$E$. Let $g\in \t C_b(E)$. Relation (\ref{4.36}) implies that there is 
a constant $b$ depending on $g$ but independent of $n$ such that 
\begin{eqnarray}\label{4.37}
\langle A_n g\, , \, A_n g\rangle_n^{1/2}\le b\, , \quad n\in {\Bbb N}. 
\end{eqnarray}
For $g\in \t C_b(E)$, set 
\begin{eqnarray*} 
\ve_n:=\left({\frac{b}{\|g\|^2}}\langle g-\beta G_{n,\beta}g\, ,
\, g-\beta G_{n,\beta}g\rangle_n^{1/2}\right)^{1/3}\, , \quad n\in 
{\Bbb N}. 
\end{eqnarray*}
Let us recall from Subsection 3.4 the definition of the set $B$ and   
condition (${\cal C}4$) where the expectation here in Section 4 is 
with respect to the initial measure $\bnu_n$, i. e., we use ${\Bbb 
E}_{\1}$ according to the notation of Subsection 3.4. By  
(\ref{4.37}) and a standard estimate on capacities (cf.~\cite{MR92}, 
V.2.6 and III.2.10), we can state 
\begin{eqnarray*} 
\, {\Bbb E}_{\1}\left(e^{-\beta\tau_{B^c}}\right)&\le&{\frac{1} 
{\ve^2_n\|g\|^2}}\cdot {\cal E}(g-\beta G_{n,\beta}g,g-\beta G_{n, 
\beta}g)+{\frac{\beta}{\ve^2_n\|g\|^2}}\cdot\langle g-\beta G_{n, 
\beta}g\, ,\, g-\beta G_{n,\beta}g \rangle_n \\ 
&=&{\frac{1}{\ve^2_n\|g\|^2}} \cdot\langle -A_ng\, ,\, g-\beta 
G_{n,\beta}g\rangle_n \\ 
&\le&{\frac{b}{\ve^2_n\|g\|^2}}\cdot\langle g-\beta G_{n,\beta} 
g \, ,\, g-\beta G_{n,\beta}g\rangle_n^{1/2} \\ 
&=&\ve_n\, . \vphantom{\frac12}
\end{eqnarray*}
Thus, we have (${\cal C}4$). 
\qed
\end{document}